\documentclass[a4paper,10pt]{article}

\usepackage{multicol}
\usepackage{ae}
\usepackage{aecompl}
\usepackage[cm]{aeguill}
\usepackage[T1]{fontenc}
\usepackage[latin1]{inputenc}
\usepackage[english]{babel}
\usepackage{fancyhdr}
\usepackage{amsmath,amssymb,amsthm}
\usepackage[all]{xy}
\usepackage{pstricks,pst-node}
\usepackage{geometry}
\newtheorem{lem}{Lemma}[section]
\newtheorem{prop}[lem]{Proposition}
\newtheorem{cor}[lem]{Corollary}
\newtheorem{definition}[lem]{Definition}
\newtheorem{rem}[lem]{Remark}
\newtheorem{ex}[lem]{Example}
\newtheorem{thm}{Theorem}[]

\def\c#1{\mathcal{#1}}
\def\square{\blacksquare}
\def\g{\Gamma}

\makeatletter
\renewcommand\section{\@startsection {section}{1}{0mm}%
                                   {-3.5ex \@plus -1ex \@minus -.2ex}%
                                   {2.3ex \@plus.2ex}%
                                   {\normalfont\large\bfseries}}
\makeatother

\begin{document}

\title{The universal cover of a monomial triangular algebra without
  multiple arrows }
\author{Patrick Le Meur
\footnote{\textit{e-mail:} patrick.lemeur@cmla.ens-cachan.fr}
\footnote{\textit{adress:} 
CMLA, ENS Cachan, CNRS, UniverSud, 61 Avenue du President Wilson, F-94230 Cachan}
}
\date{\today}

\maketitle
\abstract{
Let $A$ be a basic connected finite dimensional algebra over an
algebraically closed field $k$. Assuming that $A$ is monomial and that
the ordinary quiver $Q$ of $A$ has no oriented cycle and no multiple
arrows, we prove that $A$ admits a universal cover with group the
fundamental group of the underlying space of $Q$.
}

\section*{Introduction}
Let $A$ be a finite dimensional $k$-algebra where $k$ is an
algebraically closed field. In order to study the category $mod(A)$ of (left)
$A$-modules, one may assume that $A$ is basic and connected.
In \cite{riedtmann} (see also \cite{bongartz_gabriel}), C.~Riedtmann has introduced the covering
techniques which reduce the study of part of $mod(A)$ to
the easier one of $mod(\c C)$, where $\c C\to A$ is a Galois covering
and $\c C$ is locally bounded. These techniques
are based on the coverings of
translation quivers and their fundamental group, and therefore,
have been particularly efficient for representation-finite and
standard algebras $A$: In this case, P.~Gabriel (\cite{gabriel}) has constructed a
universal Galois covering of $A$, whose properties have led to a precise
description of the standard form of a representation-finite algebra
(\cite{bretscher_gabriel}). Unfortunately,  the above construction of
\cite{gabriel} cannot be proceeded in the representation-infinite case precisely
because the Auslander-Reiten quiver is no longer connected. 

In \cite{martinezvilla_delapena},
R.~Martinez-Villa and J.~A.~de~la~Pe\~na have constructed a Galois
covering $k\widetilde{Q}/\widetilde{I}\to A$ associated with each
presentation $kQ/I\simeq A$ (by quiver and admissible relations), for
any algebra $A$. This
Galois covering is induced by the universal cover
$(\widetilde{Q},\widetilde{I})\to (Q,I)$ with fundamental group
$\pi_1(Q,I)$ of the bound quiver $(Q,I)$, as a
generalisation of the universal cover of a translation quiver defined
in \cite{bongartz_gabriel} and \cite{riedtmann}. Like in topology, the group $\pi_1(Q,I)$ is defined
by means of a homotopy relation $\sim_I$ on the set of unoriented
paths of $Q$. When $A$ is
representation-finite and standard, this Galois covering coincides
with the one constructed by P.~Gabriel. Therefore, it is a natural
candidate for a universal cover of $A$ in the general
case. Alas, different presentations
may have non-isomorphic fundamental groups. So there may exist many
candidates for a universal cover. As an example, let
$A=kQ/<da>$, where 
$Q$ is the quiver:
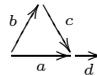
\begin{equation}
  \xymatrix@1@=10pt{&\ar@{->}[rd]^c &&\\
\ar@{->}[rr]_a \ar@{->}[ru]^b &&\ar@{->}[r]_d&}\notag
\end{equation}
Then, $\pi_1(Q,<da>)\simeq\mathbb{Z}$. On the other hand,
$A\simeq kQ/<da-dcb>$, and $\pi_1(Q,<da-dcb>)=1$. Notice that $A$ is tilted of euclidean
an type and therefore belongs to a quite well-understood class of
algebras. This illustrates the fact that except for
representation-finite algebras there are quite a few classes of algebras
for which the existence of a universal cover is known. 

In this text, we prove the existence of a universal cover for
certain monomial algebras, that is, quotients of paths algebras of
quivers by a monomial ideal (\textit{i.e.} generated by a set of
paths). More precisely, we prove the following main result.
\begin{thm}
\label{Thm1}
Let $A=kQ/I_0$, where $Q$ is a quiver without otiented cycle and without multiple
arrows, and $I_0$ is a monomial admissible ideal of $kQ$. Let $\widehat{\c C}\to
kQ/I_0$ be the Galois covering with group $\pi_1(Q)$ defined by the
presentation $kQ/I_0\simeq A$ (see
\cite{martinezvilla_delapena}). Then $\widehat{\c C}\to A$ is a
universal cover of $A$ in the following sense. For any Galois covering $\c C\to
A$ with group $G$ and with $\c C$ connected and locally bounded,
there exists a
commutative "factorisation diagram'':
\begin{equation}
\xymatrix@=8pt{
\widehat{\c C} \ar@{->}[rrd] \ar@{->}[dd] && \\
&&\c C \ar@{->}[d]\\
A \ar@{->}[rr]^{\sim} &&A
}\notag
\end{equation}
where the bottom arrow is an isomorphism of
$k$-algebras, extending the identity map on the set $Q_0$ of
vertices, and where $\widehat{\c C}\to \c C$ is a Galois covering with
group $N\vartriangleleft\pi_1(Q)$ such that there exists an
exact sequence of groups: $1\to N\to \pi_1(Q)\to G\to 1$.
\end{thm}

The author gratefully acknowledges an anonymous referee for pointing out
the following example from \cite[3.2]{geiss_delapena}. It shows
that the assumption on multiple arrows cannot be removed: Let
$A=k\left(\xymatrix@1{\cdot\ar@{->}@<2pt>[r]^x\ar@{->}@<-2pt>[r]_{x'} &
  \cdot\ar@{->}@<2pt>[r]^y\ar@{->}@<-2pt>[r]_{y'} & \cdot}\right)/I_0$, where
$I_0=<yx',y'x>$. Then, $\pi_1(Q,I_0)=\pi_1(Q)\simeq
\mathbb{Z}\star\mathbb{Z}$, and if $char(k)\neq 2$, then $A\simeq kQ/I$
where $I=<yx-y'x',yx'-y'x>$, and $\pi_1(Q,I)\simeq
\mathbb{Z}/2\mathbb{Z}$. Now, let $\c C\to A$ be the Galois covering
with group $\pi_1(Q,I)$ defined by the presentation $A\simeq kQ/I$. Then, with the notations
of Theorem~\ref{Thm1}, it is easy to show that there is no $k$-linear
functor $\widehat{\c C}\to \c C$. Thus, in this example, $A$ admits no
universal cover in the sense of Theorem~\ref{Thm1}. From the very
definition of monomial algebras, one would expect that they all have a
universal cover. The above counter-example shows that this is not
always the case. It is all the more surprising as the involved algebra
is gentle, so its representation theory its fairly well-known.\\

We now explain the strategy of the proof of Theorem~\ref{Thm1}.
Our main tool is
the quiver $\Gamma$ of the homotopy 
relations $\sim_I$ of the presentations $kQ/I\simeq A$. It was introduced in
\cite{lemeur2} to prove the existence of a universal cover for
algebras over a zero characteristic field, and whose ordinary quiver
have no double bypass. In general, if $kQ/I\simeq A$ and $kQ/J\simeq
A$, there is no simple relation between $\pi_1(Q,I)$ and
$\pi_1(Q,J)$ (and therefore between the associated Galois coverings of
$A$) unless $A$ is representation-finite (in which case $A$ is
schurian, so that $\pi_1(Q,I)=\pi_1(Q,J)$). This is the main difficulty in proving the existence
of a universal cover. Hopefully, such a relation exists when $I$ and $J$
are related by a transvection $\varphi_{\alpha,u,\tau}$, that is,
$J=\varphi_{\alpha,u,\tau}(I)$, where $(\alpha,u)$ is a bypass
(meaning that $\alpha$ is an arrow and $u$ is a path parallel to and different from
$\alpha$), $\tau\in k$ and $\varphi_{\alpha,u,\tau}\in Aut(kQ)$ is the
automorphism which maps $\alpha$ to $\alpha+\tau u$ and which fixes
any other arrow. In such a case, there is a natural quotient relation
between  $\pi_1(Q,I)$ and $\pi_1(Q,J)$.
Besides, the Galois coverings of
$A$ with groups $\pi_1(Q,J)$ and $\pi_1(Q,I)$ are the vertical arrows
of a factorisation diagram like in Theorem~\ref{Thm1}, and the
associated exact sequence of groups is given by the above quotient relation. The quiver
$\Gamma$ is then defined as follows. Its vertices are the homotopy
relations $\sim_I$ of all the presentations $kQ/I\simeq A$, and there
is an arrow $\sim\to \sim'$ if there exist presentations $kQ/I\simeq
A$ and $kQ/J\simeq A$, and a transvection $\varphi$ such that:
$\sim=\sim_I$, $\sim'=\sim_J$, $J=\varphi(I)$, and $\pi_1(Q,J)$ is a
strict quotient of $\pi_1(Q,I)$. The quiver $\Gamma$ is then
finite, connected, and it has no oriented cycle. Notice that
$\Gamma$ is reduced to a point (with no arrow) when $A$ is schurian
(and in particular when $A$ is representation-finite). But, usually,
$\Gamma$ has many vertices and many arrows.  We refer the reader
to Section~$1$ for more
details.

Roughly speaking, the existence of a universal cover is related to the
existence of a unique source in $\Gamma$. More precisely, assume
that there exists a presentation $kQ/I_0\simeq A$ (which in our case
will be the monomial presentation) such that for any other
presentation $kQ/I\simeq A$, there exists a sequence of ideals
$I_0,I_1=\varphi_1(I_0),\ldots,I_n=\varphi_n(I_{n-1})=I$, where
$\varphi_1,\ldots,\varphi_n$ are transvections defining
a path $\sim_{I_0}\to \sim_{I_1}\to\ldots\to \sim_{I_n}=\sim_I$ in $\Gamma$. Then,
the Galois covering of $A$ with group $\pi_1(Q,I_0)$
associated to the presentation  $kQ/I_0\simeq A$ is a universal cover
of $A$.
As an example, assume that $A=kQ/I_0$, where $Q$ is the quiver
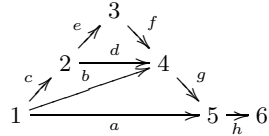
\begin{equation}
  \xymatrix@=8pt{
    &&3\ar@{->}[rd]^f&&&\\
    &2\ar@{->}[ru]^e \ar@{->}[rr]^d&&4\ar@{->}[rd]^g&&\\
    1\ar@{->}[ru]^c \ar@{->}[rrru]^b \ar@{->}[rrrr]_a&&&&5\ar@{->}[r]_h&6
  }\notag
\end{equation}
and $I_0=<ha,gb,dc>$. Then $\Gamma$ has the following
shape:
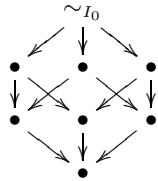
\begin{equation}
\xymatrix@=10pt{
  & \sim_{I_0}  \ar@{->}[d] \ar@{->}[ld] \ar@{->}[rd]&\\
\bullet\ar@{->}[d] \ar@{->}[rd]&
\bullet\ar@{->}[ld] \ar@{->}[rd]&
\bullet\ar@{->}[d] \ar@{->}[ld]\\ 
\bullet\ar@{->}[rd]
&\bullet\ar@{->}[d]&\bullet\ar@{->}[ld]\\ 
&\bullet&}\notag
\end{equation}
We do not specify all the vertices, but one can check that for every
$\sim_I\in\Gamma$, the group $\pi_1(Q,I)$ is free over $3-l$
generators, where $l$ is the length of any path from $\sim_{I_0}$ to
$\sim_I$ (so two presentations may have distinct homotopy relations
yet isomorphic fundamental groups). For the
needs of the proof, we construct a specific total
order on the set of bypasses of $Q$. In our example, this order is:
$(d,fe)<(b,fec)<(b,dc)<(a,gfec)<(a,gdc)<(a,gb)$. Now, let
$I=<ha+hgfec,gb+gfec,dc>$, then:
\begin{enumerate}
\item $I=\varphi_{a,gfec,1}\varphi_{b,fec,1}(I_0)$. Moreover,
  $\varphi_{a,gfec,1}\varphi_{b,fec,1}$ is the unique automorphism of
  $kQ$ which transforms $I_0$ into $I$, and which maps every arrow
  $\alpha$ to the sum of $\alpha$ and a linear combination of paths of
  length greater than $1$, none of which lying in $I_0$ (indeed:
  $gfec,fec\not\in I_0$). For that reason, we set
  $\psi_I:=\varphi_{a,gfec,1}\varphi_{b,fec,1}$.
\item The equality $\psi_I=\varphi_{a,gfec,1}\varphi_{b,fec,1}$
  expresses $\psi_I$ as a product of transvections
  $\varphi_{a,gfec,1}, \varphi_{b,fec,1}$. The associated
  sequence of bypasses is decreasing ($(a,gfec)>(b,fec)$). Actually,
  the sequence $\varphi_{a,gfec,1}, \varphi_{b,fec,1}$ is unique for
  this property.
\item We have a path $I_0\to \varphi_{b,fec,1}(I_0)\to
  \varphi_{a,gfec,1}=\varphi_{b,fec,1}(I_0)=I$ in $\Gamma$.
\end{enumerate}
In this example, it is easy to check that for any other
presentation $A\simeq kQ/J$, the ideal $J$ defines a unique
automorphism $\psi_J$ (as in 1.) which decomposes uniquely (as in 2.),
giving rise to a path from $\sim_{I_0}$ to $\sim_J$ (as in 3.).

The proof of Theorem~\ref{Thm1} mimickes the three steps we proceeded
in that example. Indeed, we prove the three following
technical points:
\begin{enumerate}
\item If $kQ/I\simeq kQ/I_0$, then there exists a unique product
  $\psi_I$ of
  transvections 
  such that $\psi_I(I_0)=I$, and such that $\psi_I$ maps every arrow
  $\alpha$ to the sum of $\alpha$ and a linear combination of paths
  of length greater than $1$,
  none of which lying in $I_0$.
\item There exists a suitable ordering on the set of bypasses such
  that if $\psi\in Aut(kQ)$ is a
  product of transvections, then $\psi$ can be written uniquely as
  $\psi=\varphi_{\alpha_n,u_n\tau_n}\ldots\varphi_{\alpha_1,u_1,\tau_1}$
  with 
  $\tau_1,\ldots,\tau_n\in k^*$ and $(\alpha_n,u_n)>\ldots>(\alpha_1,u_1)$.
\item If $kQ/I\simeq kQ/I_0$, the unique ordered sequence of
  transvections given by $1$. and $2$. yield a path in $\Gamma$
  starting at $\sim_{I_0}$ and ending at $\sim_I$. Also, this sequence
  gives rise to the factorisation diagram of Theorem~\ref{Thm1}.
\end{enumerate}

The text is therefore organised as follows. In Section~$1$ we recall all the
notions that we need to prove Theorem~\ref{Thm1}. In Section~$2$,
we prove some combinatorial facts on the paths in  a quiver. These
 lead to the order and to the decomposition of the second point
above. In Section~$3$ we prove the first point above. Finally, in
Section~$4$ we prove the last point and Theorem~\ref{Thm1}.

\section{Basic definitions}
A \textbf{$k$-category} is a category $\c C$
whose objects class $\c
C_0$ is a set, whose space of morphisms from $x$ to $y$ (denoted by $_y\c C_x$) 
is a $k$-vector space for any $x,y\in\c C_0$ and whose
composition of morphisms is $k$-bilinear. All functors between
$k$-categories will be assumed to be $k$-linear
functors. In particular, $Aut(\c C)$ will denote the group of
$k$-linear automorphism of $\c C$, and $Aut_0(\c
C)$ will denote by  for the subgroup $\{\psi\in Aut(\c C)\ |\ \psi(x)=x\ \text{for
any}\ x\in \c C_0\}$ of $Aut(\c C)$. The $k$-category $\c
C$ is called \textbf{connected} if it cannot be written as the
disjoint union of two full subcategories. An ideal $I$ of $\c C$ is the
data of subspaces $_yI_x\subseteq\ _y\c C_x$ (for any $x,y\in\c C_0$)
such that $fgh\in I$ whenever $f,g,h$ are composable morphisms
in $\c C$ such that $g\in I$.
The $k$-category $\c C$ is called \textbf{locally bounded} provided that:
$1$) for any $x\in\c C_0$, the vector spaces
  $\bigoplus\limits_{y\in\c C_0}\ _y\c C_x$ and
  $\bigoplus\limits_{y\in\c C_0}\ _x\c C_y$ are finite dimensional,
$2$) $_x\c C_x$ is a local algebra for any $x\in\c C_0$,
$3$) distinct objects are not isomorphic.
Let $A$ be a finite dimensional $k$-algebra and
let $\{e_1,\ldots,e_n\}$ be a complete set of
primitive orthogonal idempotents. Then $A$ is also a
$k$-category: $A_0:=\{e_1,\ldots,e_n\}$, $_{e_i}A_{e_i}:=e_jAe_i$ and the
product of $A$ induces the composition of morphisms.
Notice that different choices for the idempotents $e_1,\ldots,e_n$ give rise to
isomorphic $k$-categories. Also, $A$ is connected (resp. basic) as a $k$-algebra 
if and only if it is connected (resp. locally bounded) as a
$k$-category. In the sequel we shall make no distinction between a
finite dimensional $k$-algebra and its associated $k$-category.
If $\c C$ is a locally bounded $k$-category, the radical of $\c C$ is
the ideal $\c R\c C$ of $\c C$ such that: $_y\c R\c C_x$ is
the space of non-isomorphisms $x\to y$ in $\c C$, for any $x,y\in\c C_0$.
The ideal of $\c C$ generated by compositions $gf$ where $f$ and $g$
lie in $\c R\c C$ will be denoted by $\c R^2\c C$.\\

A \textbf{Galois covering
with group $G$} of $\c C$ (by $\c C'$) is a functor $F\colon\c C'\to\c C$
endowed with a group morphism $G\to Aut(\c C')$ and such that: $1$) the
induced action of $G$ on $\c C'_0$ is free, $2$) $F\circ g=F$ for any $g\in
G$, $3$) for any $k$-linear functor $F'\colon \c C'\to \c C''$ such that
$F'\circ g=F'$ for any $g\in G$, there exists a unique
$\overline{F'}\colon \c C\to \c C''$ such that $\overline{F'}\circ
F=F'$ (in other words, $F$ is a quotient of $\c C'$ by $G$ in the
category of $k$-categories). For short, the Galois covering $F$ is
called connected if $\c C'$ is
 connected and locally bounded (this implies that $\c C$ is
connected and locally bounded).
For more details on Galois coverings (in particular for the
connections with representations theory), we refer the reader to
\cite{bongartz_gabriel}.\\

\textbf{Quivers, paths, bypasses}. A quiver is a $4$-tuple
$Q=(Q_1,Q_0,s,t)$ where $Q_1$ and $Q_0$ are sets
and $s,t\colon Q_1\to Q_0$ are maps. The elements of $Q_1$ (resp. of
$Q_0$) are called the arrows (resp. the vertices) of $Q$. If
$\alpha\in Q_1$, the vertex $s(\alpha)$ (resp. $t(\alpha)$) is called
the source (resp. the target) of $\alpha$. The quiver
$Q$ is called \textbf{locally finite} if and only if any vertex is the
source (resp. the target) of finitely many arrows. 
For example, if $\c C$ is a locally bounded $k$-category, the
\textbf{ordinary quiver of $\c C$} is the locally finite quiver
$Q$ such that: $Q_0:=\c C_0$ and for
any $x,y\in \c C_0$, the number of arrows starting at $x$ and
arriving at $y$ is equal to $dim_k\ _y\c R\c C_x/\,_y\c R^2\c C_y$.
A path in $Q$ of length $n$
($n\geqslant 0$)
 with source $x\in Q_0$ (or starting at $x$) and target $y\in Q_0$ (or
 arriving at $y$) is a sequence
of arrows $\alpha_1,\ldots,\alpha_n$ such that: $x=y$ if $n=0$,
$s(\alpha_1)=x$, $s(\alpha_{i+1})=t(\alpha_i)$ for any
$i\in\{1,\ldots,n-1\}$ and $t(\alpha_n)=y$. If $n\geqslant 1$ this
path will be written $\alpha_n\ldots\alpha_1$ and called non
trivial. If $n=0$ this path will be written $e_x$ and called
stationary at $x$. The length of this path is $|u|:=n$. The mappings $s,t$ are naturally extended to paths
in $Q$. 
If $u$ and $v$ are paths, the concatenation $vu$ is defined if
and only if $t(u)=s(v)$ by the following rule:
$1$) $vu=v$ is $u$ is stationary,
$2$) $vu=u$ is $v$ is stationary,
$3$) $vu=\beta_m\ldots\beta_1\alpha_n\ldots\alpha_1$ if
  $v=\beta_m\ldots\beta_1$ and $u=\alpha_n\ldots\alpha_1$ (with
  $\alpha_i,\beta_j\in Q_1$).
Two paths in $Q$ are called \textbf{parallel} whenever they have the
same source and the same target.
 An \textbf{oriented cycle} in $Q$ is a non trivial path whose source and target are
 equal. The quiver $Q$ is said to have \textbf{multiple arrows} if and only if there
 exist in $Q$ distinct parallel arrows. If $Q$ has no oriented cycle
 and if $(\alpha,u,\beta,v)$ is a double bypass (see the introduction)
  there exists two unique paths $u_1,u_2$ such that $u=u_2\beta
  u_1$. In such a situation, the path $u_2vu_1$ will be called obtained from
 $u=u_2\beta u_1$ after replacing $\beta$ by $v$. 
\\

\textbf{Admissible presentations (see \cite[2.1]{bongartz_gabriel})}.
A quiver $Q$ defines the \textbf{path category} $kQ$ such that
$(kQ)_0=Q_0$, such that $_ykQ_x$ is the
$k$-vector space with basis the family of paths starting at $x$
and arriving at $y$, and the composition in $kQ$ is induced by the
concatenation of paths.
For short, \textbf{a normal form for $r\in\ _ykQ_x$} is an equality
$r=\sum\limits_{i=1}^n t_i u_i$ where $t_1,\ldots,t_n\in k^*$ and
$u_1,\ldots u_n$ are pairwise distinct paths in
 $Q$. With this notation, \textbf{the support} of $r$
 is the set $supp(r):=\{u_1,\ldots,u_n\}$ ($supp(0)=\emptyset$). A
\textbf{subexpression of $r$}  is a linear combination $\sum\limits_{i\in
 E}t_i u_i$ with 
 $E\subseteq\{1,\ldots,n\}$.
Later, we will need the following fact: if $r=r_1+\ldots+r_n\in\ _ykQ_x$ is such that
 $supp(r_1),\ldots,supp(r_n)$ are pairwise disjoint, then
 $r_{i_1}+\ldots+r_{i_t}$ is a subsexpression of $r$, for any indices
 $1\leqslant i_1<\ldots<i_t\leqslant n$.
 An
ideal $I$ of $kQ$ is called \textbf{admissible} provided that:
$1$) any morphism in $I$ is a linear combination of paths of length at
least $2$,
$2$) the factor category $kQ/I$ is locally bounded.
A morphism in $I$ is called a \textbf{relation} (of $I$). In particular,
\textbf{a minimal relation of $I$} (see \cite{martinezvilla_delapena})
is a non zero relation $r$ of $I$ such that $0$  and $r$ are the only
subexpressions of $r$ which are relations. With this definition, any
relation of $I$ is the sum of minimal relations with pairwise disjoint
supports. A \textbf{monomial relation} is a path $u$ lying in $I$ and $I$ is called
monomial if it is generated by a set of monomial relations.
A pair $(Q,I)$ where $Q$ is a locally finite and $I$ is an admissible
ideal of $kQ$ is called a \textbf{bound quiver}. In such a case,
$kQ/I$ is locally bounded and it is connected if and only
if $Q$ is connected (i.e. the underlying graph of $Q$ is
connected). Conversely, if $\c C$ is a locally bounded $k$-category,
then there exists an isomorphism $kQ/I\xrightarrow{\sim}\c C$ where
$(Q,I)$ is a bound quiver such that $Q$ is the ordinary quiver of $\c
C$. Such an isomorphism is called \textbf{admissible presentation of
$\c C$}. If the ideal $I$ is monomial, the admissible presentation and
$\c C$ are called monomial. Notice that $\c C$ may have different admissible
presentations.\\

\textbf{Fundamental group of a presentation (see
\cite{martinezvilla_delapena})}. Let $(Q,I)$ be a bound quiver and let
$x_0\in Q_0$.
For every arrow $x\xrightarrow{a}y\in Q_1$
we define its formal inverse $a^{-1}$ with source
$s(a^{-1})=y$ and target $t(a^{-1})=x$. This defines a new quiver
$\overline{Q}=(Q_0,Q_1\cup\{a^{-1}\ |\ a\in Q_1\},s,t)$. With these notations, \textbf{a
walk in $Q$} is a path in $\overline{Q}$. The concatenation of walks
in $Q$ is by definition the concatenation of paths in $\overline{Q}$. The
\textbf{homotopy relation} of $(Q,I)$ is the equivalence relation on
the set of walks in $Q$, denoted by $\sim_I$ and generated by the
following properties:
\begin{enumerate}
\item $\alpha\alpha^{-1}\sim_I e_y$ and $\alpha^{-1}\alpha\sim_I e_x$
for any arrow $x\xrightarrow{\alpha}y$ in $Q$,
\item $u\sim_I v$ for any $u,v\in supp(r)$ where $r$ is a minimal
relation of $I$,
\item $wvu\sim_I wv'u$ for any walks $w,v,v',u$ such that $v\sim_I v'$
and such that the concatenations $wvu$ and $wv'u$ are well-defined
(i.e. $\sim_I$ is compatible with the concatenation).
\end{enumerate}
The $\sim_I$-equivalence class of a walk $\gamma$ will be denoted by $[\gamma]_I$.
Let $\pi_1(Q,I,x_0)$ be the set of equivalence classes of walks in $Q$
with source and target equal to $x_0$. The concatenation of walks
endows this set with a group structure (with unit $e_{x_0}$) and this group is called the
\textbf{fundamental group of $(Q,I)$ at $x_0$}. If $Q$ is connected, the
isomorphism class of this group does not depend on $x_0\in
Q_0$ and $\pi_1(Q,I,x_0)$ is denoted by $\pi_1(Q,I)$. If
$\c C$ is a connected locally bounded $k$-category and if $kQ/I\simeq
\c C$ is an admissible presentation, the fundamental group
$\pi_1(Q,I)$ is called the fundamental group of this presentation.\\

\textbf{Dilatations, transvections (see \cite{lemeur2})}. Let $Q$ be a
quiver. A \textbf{dilatation} of $kQ$ is an automorphism $D\in Aut_0(kQ)$ such that
$D(\alpha)\in k^*\alpha$ for any $\alpha\in Q_1$. The dilatations of $kQ$ form a
subgroup $\c D$ of $Aut_0(kQ)$. Let $(\alpha,u)$ be a
bypass in $Q$ and let $\tau\in k$. This defines
$\varphi_{\alpha,u,\tau}\in Aut_0(kQ)$ as follows:
$\varphi_{\alpha,u,\tau}(\alpha)=\alpha+\tau u$ and
$\varphi_{\alpha,u,\tau}(\beta)=\beta$ for any arrow $\beta\neq
\alpha$. The automorphism $\varphi_{\alpha,u,\tau}$ is called a
\textbf{transvection}. The composition of transvections is ruled as follows.
Let $\varphi_{\alpha,u,\tau}$ and $\varphi_{\alpha,u,\tau'}$, then
$\varphi_{\alpha,u,\tau}\varphi_{\alpha,u,\tau'}=\varphi_{\alpha,u,\tau+\tau'}\
\ \text{and}\ \ \varphi_{\alpha,u,\tau}^{-1}=\varphi_{\alpha,u,-\tau}$.
 If $(\alpha,u,\beta,v)$ and $(\beta,v,\alpha,u)$ are not a double bypasses, then
$\varphi_{\alpha,u,\tau}\varphi_{\beta,v,\nu} =
\varphi_{\beta,v,\nu}\varphi_{\alpha,u,\tau},$.  If $(\alpha,u,\beta,v)$ is a double bypass and if $Q$ has no oriented cycle, then
$\varphi_{\beta,v,\nu}\varphi_{\alpha,u,\tau} =
\varphi_{\alpha,u,\tau}\varphi_{\alpha,w,\tau\nu}\varphi_{\beta,v,\nu}
$, where $w$ is the path obtained from $u$ after replacing $\beta$ by
$v$.
The subgroup of $Aut_0(kQ)$ generated by all the
transvections is denoted by $\c T$. The dilatations and the transvections are useful to
compare admissible presentations of an algebra because of the following proposition:
\begin{prop}
\label{prop1.1}
(see \cite[Prop. 2.1, Prop. 2.2]{lemeur2})
Let $kQ/I\simeq A$ and $kQ/J\simeq A$ be admissible presentations of
the basic finite dimensional algebra $A$. If $Q$ has no oriented
cycle, then there exists $\psi\in Aut_0(kQ)$ such that $\psi(I)=J$.
Moreover, 
$\c T$ is a normal subgroup of $Aut_0(kQ)$ and $Aut_0(kQ)=\c T\c D=\c
D\c T$.
\end{prop}

The dilatations and the transvections were introduced because they
allow comparisons between the fundamental groups of presentations of
the same locally bounded $k$-category. Notice that if $I,J$ are
admissible ideals of $kQ$ such that $\gamma\sim_I\gamma'\Rightarrow
\gamma\sim_J\gamma'$ for any walks $\gamma,\gamma'$, then the identity
map on the set of walks induces a surjective group morphism
$\pi_1(Q,I)\twoheadrightarrow\pi_1(Q,J)$. In particular, if $\sim_I$
and $\sim_J$ coincide, then $\pi_1(Q,I)=\pi_1(Q,J)$.
\begin{prop}
\label{prop1.2}
(see \cite[Prop. 2.5]{lemeur2})
Let $I$ be an
admissible ideal of $kQ$ (with $Q$ without oriented
cycle), let $\varphi\in
Aut_0(kQ)$  and set $J=\varphi(I)$. If $\varphi$ is a dilatation,
then $\sim_I$ and $\sim_J$ coincide. If
$\varphi=\varphi_{\alpha,u,\tau}$ is a transvection, then:
\begin{enumerate}
\item if $\alpha\sim_Iu$ and $\alpha\sim_Ju$ then $\sim_I$ and
$\sim_J$ coincide.
\item if $\alpha\not\sim_Iu$ and $\alpha\sim_J u$ then $\sim_J$
is generated by $\sim_I$ and  $\alpha\sim_Ju$.
\item if $\alpha\not\sim_Iu$ and $\alpha\not\sim_Ju$ then $I=J$
and $\sim_I$ and $\sim_J$ coincide.
\end{enumerate}
If there exists a transvection $\varphi$ such that $\varphi(I)=J$ and
such the second point above occurs, then 
\textbf{$\sim_J$ is called a direct successor of $\sim_I$}.
\end{prop}
Here the expression ``$\sim_J$ is generated by $\sim_I$ and
$\alpha\sim_J u$'' means that $\sim_J$ is the equivalence relation
on the set of walks in $Q$, compatible with the concatenation and
generated by the two following properties:
$1$) $\gamma\sim_I\gamma'\Rightarrow \gamma\sim_J\gamma'$,
$2$) $\alpha\sim_J u$. Following \cite[Def. 2.7]{lemeur2}, if $A$ is a
basic connected finite dimensional algebra with ordinary quiver $Q$
without oriented cycle, we define
\textbf{the quiver $\Gamma$ of the homotopy relations of $A$} to be
the quiver such that $\Gamma_0=\{\sim_I\ |\ kQ/I\simeq A\}$ and such
that there exists arrow $\sim_I\to \sim_J$ if and only if $\sim_J$ is
a direct successor of $\sim_I$. Recall (\cite[Rem. 5,
Prop. 2.8]{lemeur2}) that $\Gamma$ is finite, connected, without
oriented cycle and such that for any (oriented) path with source
$\sim_I$ and target $\sim_J$, the identity map on the set of walks in $Q$
induces a surjective group morphism $\pi_1(Q,I)\twoheadrightarrow \pi_1(Q,J)$.\\ 

\textbf{Gröbner bases} Let $E$ be a $k$-vector space with an ordered
basis $(e_1,\ldots,e_n)$, let $(e_1^*,\ldots,e_n^*)$ be the associated
dual basis of $E^*$, and let $F$ be a subspace of $E$. A Gröbner basis (see
\cite{adams_loustaunau} for the ususal definition) of $F$ is a basis
$(r_1,\ldots,r_d)$ such that:
\begin{enumerate}
\item $r_j\in e_{i_j}+Span(e_l\ ;\ l<i_j)$ for some $i_j$, for any
$j\in\{1,\ldots,r\}$,
\item $i_1<i_2<\ldots<i_r$,
\item $e_{i_j}^*(r_{j'})=0$ for any $j\neq j'$.
\end{enumerate}
It is well known that $F$ admits a unique Gröbner basis. Also, $r\in
F$ if and only if:
$r=\sum\limits_{j=1}^de_{i_j}^*(r)r_j$.
In the sequel, we will use this notion in the following setting: $E$
is the vector space with basis (for some order to be defined) the
family of non trivial paths in a
finite quiver $Q$ without oriented cycles and $F$ is the underlying
subspace of $E$ associated to an admissible ideal $I$ of $kQ$. Notice
that in this setting, the Gröbner basis of $F$ is made of minimal
relations of $I$. Also, if $r\in E$ and if $u$ is a non trivial path,
then: $u\in supp(r)\Leftrightarrow u^*(r)\neq 0$.\\

Until the end of the text, $Q$ will denote a finite quiver without
oriented cycle and without multiple arrows.

\section{Combinatorics on the paths in a quiver}
Recall from the previous section that if $(\alpha,u,\beta,v)$ is a
double bypass and if $\tau,\nu$ are scalars, then
$\psi:=\varphi_{\beta,v,\nu}\varphi_{\alpha,u,\tau}$ is equal to
$\varphi_{\alpha,u,\tau}\varphi_{\alpha,w,\tau\nu}\varphi_{\beta,v,\nu}$ 
where $w$ is the path obtained from $u$ by replacing $\beta$ by
$v$. Remark that
$\varphi_{\beta,v,\nu}\varphi_{\alpha,u,\tau}(\alpha)=\alpha+\tau
u+\tau\nu w$. Hence, the paths ($u$ and $w$) appearing in
$\varphi_{\beta,v,\nu}\varphi_{\alpha,u,\tau}(\alpha)-\alpha$ are axactely
those paths $\theta$ such that $(\alpha,\theta)$ is a bypass appearing
in one of the transvections of the product
$\varphi_{\alpha,u,\tau}\varphi_{\alpha,w,\tau\nu}\varphi_{\beta,v,\nu}$.
Moreover, the scalars ($\tau$ and $\tau\nu$) appearing with these
paths are exactely the scalars of the corresponding transvections in
this product. So, the computation
of $\varphi_{\beta,v,\nu}\varphi_{\alpha,u,\tau}(\alpha)$ can be done
just by looking for the occurences of $\alpha$ in the product
$\varphi_{\alpha,u,\tau}\varphi_{\alpha,w,\tau\nu}\varphi_{\beta,v,\nu}$.
From this point of view, the decomposition
$\psi=\varphi_{\alpha,u,\tau}\varphi_{\alpha,w,\tau\nu}\varphi_{\beta,v,\nu}$
is more useful than the decomposition $\psi=\varphi_{\beta,v,\nu}\varphi_{\alpha,u,\tau}$.
The aim of this section is to show that this phenomenon is a general one. In
this purpose the useful notion of derivation of a path and a total
order on the set of bypasses will be introduced.

\subsection{Derivation of paths}
\begin{definition}
\label{def2.1}
Let $u=\alpha_n\ldots\alpha_1$ and $v$ be paths in $Q$. Then $v$ is
called derived of $u$ (of order $t$) if there exist indices
$1\leqslant i_1<\ldots<i_t\leqslant n$ and bypasses
$(\alpha_{i_1},v_1),\ldots,(\alpha_{i_t},v_t)$ such that $v$ is
obtained from $u$ by replacing $\alpha_{i_l}$ by $v_l$, for each $l$:
\begin{equation}
v=\alpha_n\ldots\alpha_{i_t+1}v_t\alpha_{i_t-1}\ldots\alpha_{i_l+1}v_l\alpha_{i_l-1}\ldots
\alpha_{i_1+1}v_1\alpha_{i_1-1}\ldots\alpha_1\notag
\end{equation} 
\end{definition}
\begin{rem}
\label{rem2.2}
If $\alpha\in Q_1$, then $u$ is derived of $\alpha$ if and only if
$(\alpha,u)$ is a bypass.
\end{rem}

With the above definition, the following lemma is easily
verified using the fact that $Q$ has no
multiple arrows and no oriented cycle.
\begin{lem}
\label{lem2.3}
\begin{enumerate}
\item If $v$ is derived of $u$ with both orders $t$ and $t'$, then
$t$=$t'$.
\item If $v$ is derived of $u$ of order $t$ then there exists a
  sequence of paths $u_0=u,u_1,\ldots,u_t=v$ such that $u_i$ is derived of
  $u_{i-1}$ of order $1$ for any $i$.
\item If $v$ is derived of $u$ of order $t$, then $|v|\geqslant
  |u|+t$.
\item If $v$ is derived of $u$ of order $t$ and if $w$ is derived of
  $v$ of order $t'$, then $w$ is derived of $u$ of order at least $t$.
\item Let $u,v,w$ be paths verifying:
\begin{itemize}
\item $v$ is derived of $u$,
\item $w$ is derived of $v$,
\item $w$ is derived of $u$ of order $1$,
\end{itemize}
then we have:
\begin{equation}
u=u_2\alpha u_1,\ \ v=u_2\theta u_1,\ \ w=u_2\theta' u_1\notag
\end{equation}
where $u_1,u_2$ are paths, $(\alpha,\theta)$ is a bypass and
$\theta'$ is derived of $\theta$.
\item If $v$ (resp. $v'$) is derived of $u$ (resp. of $u'$) of order
  $t$ (resp. $t'$), then $v'v$ is derived of $u'u$ of order $t'+t$,
  whenever these compositions of paths are well defined.
\end{enumerate}
\end{lem}

The following example shows that, in the preceding lemma, the inequality in the $4$-th point may be an equality.
\begin{ex}
\label{ex2.4}
Let $(\alpha,u,\beta,v)$ be a double bypass. Let $u_1,u_2$ be the paths
such that $u=u_2\beta u_1$. Then $u$ is derived of $\alpha$ of order
$1$, $w:=u_2v u_1$ is derived of $u$ of order $1$ and $w$ is
derived of $u$ of order $1$.
\end{ex}

\subsection{Order between paths, order between bypasses}

Now, we construct a total order on the set of non
trivial paths in $Q$. This construction is a particular case of the
one introduced in \cite{farkas_feustel_green}. Also it depends on an
arbitrary order $\vartriangleleft$ on $Q_1$. We assume that this order
$\vartriangleleft$ is fixed for this subsection. We shall write
$\vartriangleleft$ for the lexicographical order induced by
$\vartriangleleft$ on the set of nontrivial paths in $Q$. For details
on the correctness of the following definition we refer the reader to
\cite{farkas_feustel_green}.
\begin{definition}
\label{def2.5}
For $\alpha\in Q_1$, set:
\begin{equation}
W(\alpha)=Card(B(\alpha))\ \ \text{where}\ \ B(\alpha)=\{(\alpha,u)\
|\ (\alpha,u)\ \text{is a bypass in $Q$}\}\notag
\end{equation}
For $u=\alpha_n\ldots\alpha_1$ a path in $Q$ (with $\alpha_i\in Q_1$),
let us set:
\begin{equation}
W(u)=W(\alpha_n)+\ldots+W(\alpha_1)\notag
\end{equation}
These data define a total order $<$ on the set of non trivial paths in
$Q$ as follows:
\begin{equation}
u<v\Leftrightarrow \left\{
\begin{array}{rl}
&W(u) <W(v)\\
or&\\
&W(u)=W(v)\ and\ u\vartriangleleft v
\end{array}\right.\notag
\end{equation}
We shall write $<$ for the lexicographical order induced by $<$
on the set of couples of paths.
\end{definition}

\begin{rem}
\label{rem2.6}
If
$u$ and $v$ are (non trivial) paths such that $vu$ is well defined,
then $W(vu)=W(u)+W(v)$.
\end{rem}

\begin{ex}
  Let $Q$ be the following quiver without oriented cycle and without
  multiple arrows:
  
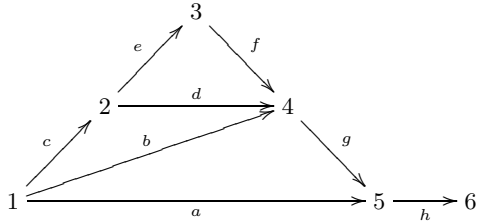
\begin{equation}
    \xymatrix{
      &&3\ar@{->}[rd]^f&&&\\
      &2\ar@{->}[ru]^e \ar@{->}[rr]^d&&4\ar@{->}[rd]^g&&\\
      1\ar@{->}[ru]^c \ar@{->}[rrru]^b \ar@{->}[rrrr]_a&&&&5\ar@{->}[r]_h&6
    }\notag
  \end{equation}
and let $\vartriangleleft$ be any total order on $Q_1$. Then, $B(a)=\{(a,bg),\,(a,gdc),\,(a,gfe)\}$,
$B(b)=\{(b,cd),\,(b,fec)\}$, $B(d)=\{(d,fe)\}$ and $B(x)=\emptyset$
for $x\in Q_1\backslash\{a,b,d\}$. In particular, the paths with
source $1$ and target $5$ are ordered as follows:
\begin{equation}
  gfec<gdc<gb<a\notag
\end{equation}
\end{ex}

\begin{lem}
\label{lem2.7}
\begin{enumerate}
\item If $u,v,u',v'$ are paths such that $v<u$
  and $v'<u'$ then $v'v<u'u$ whenever these compositions are well
  defined.
\item If $(\alpha,u)$ is a bypass, then $W(u)<W(\alpha)$. So
  $u<\alpha$.
\item If $v$ is derived of $u$, then $v<u$.
\item If $(\alpha,u,\beta,v)$ is a double bypass and if $w$ is the
  path obtained from $u$ after replacing $\beta$ by $v$, then:
\begin{equation}
(\beta,v)<(\alpha,w)<(\alpha,u)\notag
\end{equation}
\end{enumerate}
\end{lem}
\noindent{\textbf{Proof:}} $1)$ is a direct consequence of
Definition~\ref{def2.5} and Remark~\ref{rem2.6}.

$2)$ Let us write $u=a_n\ldots a_1$ with $a_i\in Q_1$ for each
$i$ (hence $a_i\neq a_j$ if $i\neq j$ because $Q$ has no oriented
cycle). Therefore:
\begin{enumerate}
\item[.] $B(a_1),\ldots,B(a_n)$ are pairwise disjoint,
\item[.] $W(u)=W(a_1)+\ldots+W(a_n)$
\end{enumerate}
Notice that if $(a_i,v)\in B(a_i)$, then
 $(\alpha,a_n\ldots a_{i+1}v a_{i-1}\ldots
a_1)\in B(\alpha)$. Thus, we have a well defined mapping:
\begin{equation}
\begin{array}{crcl}
\theta\colon & B(a_1)\sqcup\ldots\sqcup B(a_n) & \longrightarrow &
B(\alpha)\\
& (a_i,v) & \longmapsto &(\alpha,a_n\ldots a_{i+1}v a_{i-1}\ldots a_1)
\end{array}\notag
\end{equation}
This mapping is one-to-one, indeed:
\begin{enumerate}
\item[.] if $\theta(a_i,v)=\theta(a_i,v')$ with $(a_i,v),(a_i,v')\in
  B(a_i)$ then:
\begin{equation}
a_n\ldots a_{i+1}va_{i-1}\ldots a_1=a_n\ldots
  a_{i+1}v'a_{i-1}\ldots a_1\notag
\end{equation}
 and therefore $(a_i,v)=(a_i,v')$,
\item[.] if $\theta(a_i,v)=\theta(a_j,v')$ with $(a_i,v)\in B(a_i)$,
  $(a_j,v')\in B(a_j)$ and $j<i$, then:
\begin{equation}
a_n\ldots a_{i+1}va_{i-1}\ldots a_1=a_n\ldots a_{j+1}v'a_{j-1}\ldots
a_1\notag
\end{equation}
So:
\begin{equation}
va_{i-1}\ldots a_1=a_i\ldots a_{j+1}v'a_{j-1}\ldots a_1\notag
\end{equation}
Since $v$ and $a_i$ are parallel and since $Q$ has no
oriented cycle, we infer that $v=a_i$ which is impossible because
$(a_i,v)\in B(a_i)$.
\end{enumerate}
On the other hand, $\theta$ is not onto. Indeed, if there exists
$(a_i,v)\in B(a_i)$ verifying $\theta(a_i,v)=(\alpha,u)$, then:
\begin{equation}
a_n\ldots a_1=u=a_n\ldots a_{i+1}v a_{i-1}\ldots a_1\notag
\end{equation}
which implies $a_i=v$, a contradiction. Since $\theta$ is one-to-one
and not onto, we deduce that:
\begin{equation}
W(\alpha)=Card(B(\alpha))>Card(B(a_1)\sqcup \ldots\sqcup B(a_n))=W(u)\notag
\end{equation}
This proves that $W(u)<W(\alpha)$ and that $u<\alpha$.

$3)$ is a direct consequence of $1)$ and of $2)$.

$4)$ Let us write $u=u_2\beta u_1$ (with $u_1,u_2$ paths) so that
$w=u_2vu_1$. From $2)$, we have:
\begin{equation}
W(\alpha)>W(u)=W(u_1)+W(\beta)+W(u_2)\geqslant W(\beta)\notag
\end{equation}
So $\beta<\alpha$ and therefore
$(\beta,v)<(\alpha,w)$.
Using $2)$ again, we also have:
\begin{equation}
W(w)=W(u_2)+W(v)+W(u_1)<W(u_2)+W(\beta)+W(u_1)=W(u)\notag
\end{equation}
So $w<u$ and therefore
$(\alpha,w)<(\alpha,u)$
\hfill$\square$\\

Unless otherwise specified, $<$ will always denote an order on the set
of paths as in Definition~\ref{def2.5}.

\subsection{Image of a path by a product of transvections}
In this paragraph, we apply the previous constructions to find an easy
way to compute $\psi(u)$ when $\psi\in\c T$ and $u$ is a path in
$Q$. We begin with the following lemma on the description of
$\psi(\alpha)$ when $\psi\in\c T$ and $\alpha\in Q_1$. Recall that $Q$ has no
multiple arrows and no oriented cycle.
\begin{lem}
\label{lem2.8}
Let $\psi\in\c T$ and let $\alpha\in Q_1$. Then $\psi(\alpha)-\alpha$
is a linear combination of paths parallel to $\alpha$ and of length
greater than or equal to $2$. In particular, $\alpha\in
supp(\psi(\alpha))$ and $\alpha^*(\psi(\alpha))=1$.
\end{lem}
\noindent{\textbf{Proof:}}
The conclusion is immediate if $\psi$ is a transvection because $Q$
has no multiple arrows. The conclusion in the general case is obtained
using an easy induction on the number of transvections whose product
equal $\psi$.\hfill$\square$\\

The preceding lemma gives the following description of $\psi(u)$ when
$\psi\in\c T$ and $u$ is a path. We omit the proof which is immediate
thanks to Lemma~\ref{lem2.8} and to point $6)$ of Lemma~\ref{lem2.3}.
\begin{prop}
\label{prop2.9}
Let $\psi\in\c T$ and let $u=\alpha_n\ldots\alpha_1$ be a path in $Q$
(with $a_i\in Q_1$ for any $i$). For each $i$, let:
\begin{equation}
\psi(\alpha_i)=\alpha_i+\sum\limits_{j=1}^{m_i}\lambda_{i,j}u_{i_j}\notag
\end{equation}
be a normal form for $\psi(\alpha_i)$. Then $supp(\psi(u))$ is the set
of the paths in $Q$ described as follows. Let $r\in \{0,\ldots,n\}$.
  Let $1\leqslant i_1<\ldots<i_r\leqslant n$ be indices. For each
  $l\in\{1,\ldots,r\}$, let $j_l\in\{1,\ldots, m_{i_l}\}$. Then the
    following path obtained from $u$ after replacing $\alpha_{i_l}$ by
    $u_{j_l}$ for each $l$ belongs to $supp(\psi(u))$:
\begin{equation}
\alpha_n\ldots\alpha_{i_r+1}u_{j_r}\alpha_{i_r-1}
\ldots\alpha_{i_l+1}u_{j_l}\alpha_{i_l-1}\ldots
\alpha_{i_1+1}u_{j_1}\alpha_{i_1-1}\ldots\alpha_1\notag
\end{equation}
Moreover, this path appears in $\psi(u)$ with coefficient:
\begin{equation}
\lambda_{i_1,j_1}\ldots \lambda_{i_r,j_r}\notag
\end{equation}
As a consequence, $\psi(u)-u$ is a linear combination of paths derived
of $u$.
\end{prop}

\begin{ex}
The previous proposition does not hold if $Q$ has multiple arrows. For
example, if $Q$ is the Kronecker quiver
$\xymatrix{1\ar@/^/@{->}[r]^a\ar@/_/@{->}[r]_b&2}$ and if
  $\psi=\varphi_{a,b,1}\varphi_{b,a,-1}\varphi_{a,b,1}$, then
  $\psi(a)=b$ and $\psi(b)=-a$.  
\end{ex}

\begin{rem}
\label{rem2.10}
If $(\alpha,u)$ is a bypass and if $v\in supp(\psi(u)-u)$, then
$(\alpha,v)$ is also a bypass and $(\alpha,v)<(\alpha,u)$.
\end{rem}

Now we are able to state the main result of this paragraph. It describes
$\psi(\alpha)$ $(\alpha\in Q_1)$ using a particular writing of $\psi$
as a product of transvections. Notice that the following proposition
formalises the phenomenon observed at the begining of the section.
\begin{prop}
\label{prop2.11}
Let $(\alpha_1,u_1)<\ldots<(\alpha_n,u_n)$ be an increasing sequence of
bypasses, let
$\tau_1,\ldots,\tau_n\in k^*$ and set
$\psi=\varphi_{\alpha_n,u_n,\tau_n}\ldots\varphi_{\alpha_1,u_1,\tau_1}$.
For any $\alpha\in Q_1$, there is a normal form for $\psi(\alpha)$:
\begin{equation}
\psi(\alpha)=\alpha+\sum_{i\ such\ that\ \alpha=\alpha_i}\tau_iu_i\notag
\end{equation}
\end{prop}
\noindent{\textbf{Proof:}} Let us prove that the conclusion of the
proposition is true using an induction on $n\geqslant 1$. By
definition of a transvection, the proposition holds of $n=1$. Assume
that $n\geqslant 2$ and that the conclusion of the proposition holds
if we replace
$\psi=\varphi_{\alpha_n,u_n,\tau_n}\ldots\varphi_{\alpha_1,u_1,\tau_1}$
by
$\varphi_{\alpha_{n-1},u_{n-1},\tau_{n-1}}\ldots\varphi_{\alpha_1,u_1,\tau_1}$.
Therefore, for $\alpha\in Q_1$, we have a normal form:
\begin{equation}
\varphi_{\alpha_{n-1},u_{n-1},\tau_{n-1}}\ldots\varphi_{\alpha_1,u_1,\tau_1}(\alpha)=
\alpha +\sum\limits_{i\leqslant n-1,\ \alpha=\alpha_i}\tau_iu_i\notag
\end{equation}
So:
\begin{equation}
\psi(\alpha) = \varphi_{\alpha_n,u_n,\tau_n}(\alpha)+
\sum\limits_{i\leqslant n-1,\
  \alpha=\alpha_i}\tau_i\varphi_{\alpha_n,u_n,\tau_n}(u_i)\tag{$i$}
\end{equation}
Let $i\in\{1,\ldots,n-1\}$. Thanks to Lemma~\ref{lem2.7}, the
  inequality $(\alpha_i,u_i)<(\alpha_n,u_n)$ implies that
  $(\alpha_i,u_i,\alpha_n,u_n)$ is not a double bypass. Thus,
  $\alpha_n$ does not appear in the path $u_i$. This proves that:
\begin{equation}
(\forall i\in\{1,\ldots,n-1\})\ \
\varphi_{\alpha_n,u_n,\tau_n}(u_i)=u_i\tag{$ii$}
\end{equation}
The definition of $\varphi_{\alpha_n,u_n,\tau_n}$, together with $(i)$
and $(ii)$, imply the equality:
\begin{equation}
\psi(\alpha) = \alpha+
\sum\limits_{
  \alpha=\alpha_i}\tau_iu_i\tag{$iii$}
\end{equation}
It only remains to prove that the equality $(iii)$ is a normal
form. Remark that all the scalars which appear in the right-hand side
of $(iii)$ are non zero. Moreover, if $i\in\{1,\ldots,n\}$ verifies
$\alpha=\alpha_i$, then $\alpha\neq u_i$, because $(\alpha,u_i)$ is a
bypass. Finally, if $1\leqslant i<j\leqslant n$ verify
$\alpha=\alpha_i=\alpha_j$, then $(\alpha,u_i)=(\alpha_i,u_i)<(\alpha_j,u_j)=(\alpha,u_j)$ so
$u_i\neq u_j$. Therefore, $(iii)$ is a normal form for
$\psi(\alpha)$.\hfill$\square$\\

When $\psi\in \c T$ is like in Proposition~\ref{prop2.11}, we
shall say that $\psi$ is written as a decreasing product of
transvections. Later we will prove that any $\psi\in\c T$ can be
written uniquely as a decreasing product of transvections.
The description in Proposition~\ref{prop2.11} will be particularly
useful in the sequel.
We end this paragraph with two propositions concerning the description
of $\psi(r)$ when $\psi\in\c T$ and $r$ is a linear combination of
paths. The following proposition gives conditions for $\psi^{-1}(r')$
to be a subexpression of $r$ when $r'$ is a subexpression of
$\psi(r)$.
\begin{prop}
\label{prop2.12}
Let $\psi\in \c T$, let $r\in\ _ykQ_x$ and let $r'$ be a subexpression
of $\psi(r)$. Let $\simeq$ be the equivalence relation on the set of
paths in $Q$ generated by:
\begin{equation}
v\in supp(\psi(u))\Rightarrow u\simeq v\notag
\end{equation}
Assume that for any $u,v\in supp(\psi(r))$ verifying $u\simeq v$ we
have:
\begin{equation}
u\in supp(r')\Leftrightarrow v\in supp(r')\notag
\end{equation}
Then $\psi^{-1}(r')$ is a subexpression of $r$.
\end{prop}
\noindent{\textbf{Proof:}} Let $\simeq'$ be the trace of $\simeq$ on
$supp(r)$ and let us write $supp(r)=c_1\sqcup \ldots \sqcup c_n$ as a disjoint union of its
$\simeq'$-classes. This partition of $supp(r)$ defines a decomposition of $r=r_1+\ldots+r_n$
where $r_i$ is the subexpression of $r$ verifying $supp(r_i)=c_i$. For
each $i$, let us fix a normal form:
\begin{equation}
r_i=\sum\limits_{j=1}^{n_i}t_{i,j}u_{i,j}\notag
\end{equation}
so that we have the following normal form for $r$:
\begin{equation}
r=\sum\limits_{i=1}^n\sum\limits_{j=1}^{n_i}t_{i,j}u_{i,j}\notag
\end{equation}
Let us set $r_i':=\psi(r_i)$. In order to prove that $\psi^{-1}(r')$
is a subexpression of $r$, we will prove that there exist indices
$1\leqslant i_1<\ldots<i_t\leqslant n$ verifying
$r'=r_{i_1}'+\ldots+r_{i_t}'$ (so that
$\psi^{-1}(r')=r_{i_1}'+\ldots+r_{i_t}'$). In this purpose, we will
successively prove the following
facts:
\begin{enumerate}
\item[$1)$] $u,v\in supp(r_i')\Rightarrow u\simeq v$, for any $i$,
\item[$2)$] $supp(r_1'),\ldots,supp(r_n')$ are pairwise disjoint,
\item[$3)$] for each $i$, $r_i'$ is a subexpression of $\psi(r)$,
\item[$4)$] if $i\in\{1,\ldots,n\}$ verifies $supp(r')\cap
  supp(r_i')\neq\emptyset$, then $supp(r_i')\subseteq supp(r')$,
\end{enumerate}
$1)$ Let $i\in\{1,\ldots,n\}$ and let $u,v\in supp(r_i')$. So there
exist $u',v'\in supp(r_i)$ such that
$u\in supp(\psi(u'))$ and $v\in supp(\psi(v'))$.
By definition of $\simeq$ and of $r_i$, we deduce that:
\begin{equation}
u,v\in supp(r_i')\Rightarrow u\simeq v\tag{$i$}
\end{equation}
$2)$ Let $i,j\in \{1,\ldots,n\}$ be such that there exists $v\in
supp(r_i')\cap supp(r_j')$. So there exist $u\in supp(r_i)$ and $u'\in
supp(r_j)$ such that
$v\in supp(\psi(u))$ and $v\in supp(\psi(u'))$.
This implies that
$u\simeq v\simeq u'$.
Since $u\in c_i=supp(r_i)$ and $u'\in c_j=supp(r_j)$, we deduce that
$c_i=c_j$ and therefore $i=j$. So:
\begin{equation}
i\neq j\Rightarrow supp(r_i')\cap supp(r_j')=\emptyset\tag{$ii$}
\end{equation}
$3)$ We have $\psi(r)=r_1'+\ldots +r_n'$ so $(ii)$ implies that:
\begin{equation}
\text{$r_i'$ is a subexpression of $\psi(r)$ for any $i$}\tag{$iii$}
\end{equation}
$4)$ Let $i\in\{1,\ldots,n\}$ and assume that there exists $u\in
supp(r_i')\cap supp(r')$. If $v\in supp(r_i')$ then $u\simeq v$ thanks
to $(i)$. So, by assumption on $r'$, we have $v\in supp(r')$. This
proves that:
\begin{equation}
supp(r_i')\cap supp(r')\neq \emptyset\Rightarrow supp(r_i')\subseteq
supp(r')\tag{$iv$}
\end{equation}
Now, we can prove that $\psi^{-1}(r')$ is a subexpression of
$r$. Thanks to $(iii)$, the elements $r',r_1',\ldots,r_n'$ are
subexpressions of $\psi(r)$. So $(iv)$ and the equality
$\psi(r)=r_1'+\ldots +r_n'$ imply that there exist indices $1\leqslant
i_1<\ldots < i_t\leqslant n$ such that
$r'=r_{i_1}'+\ldots+r_{i_t}'$.
So
$\psi^{-1}(r')=r_{i_1}+\ldots+r_{i_n}$.
This proves that $\psi^{-1}(r')$ is a subexpression of
$r$.\hfill$\square$\\

The last proposition of this subsection gives a sufficient condition
on $u\in supp(r)$ to verify $u\in supp(\psi(r))$.
\begin{prop}
\label{prop2.13}
Let $\psi\in\c T$, let $r\in\ _ykQ_x$ and let $u\in supp(r)$. Then, at
least one of the two following facts is verified:
\begin{enumerate}
\item $u\in supp(\psi(r))$,
\item there exists $v\in supp(r)$ such that $u\neq v$ and such that
  $u\in supp(\psi(v))$.
\end{enumerate}
As a consequence, if $u$ is not derived of $v$ for any $v\in supp(r)$,
then:
\begin{equation}
u\in supp(\psi(r))\ \ \text{and}\ \ u^*(\psi(r))=u^*(r)\notag
\end{equation}
\end{prop}
\noindent{\textbf{Proof:}} Let us fix a normal form 
$r=\sum\limits_{i=1}^nt_iu_i$
where we may assume that $u=u_1$. Let us assume that $u\not\in
supp(\psi(r))$, i.e. $u^*(\psi(r))=0$. Recall from
Proposition~\ref{prop2.9} that $u^*(\psi(u))=1$, so:
\begin{equation}
0=u^*(\psi(r))=t_1+\sum\limits_{i=2}^nt_iu^*(\psi(u_i))\tag{$i$}
\end{equation}
Therefore, there exists $i_0\in\{2,\ldots,n\}$
such that $u^*(\psi(u_{i_0}))\neq 0$. So:
\begin{equation}
u_{i_0}\in supp(r),\ \ u_{i_0}\neq u_1=u\ \ \text{and}\ \
u_1^*(\psi(u_{i_0}))\neq 0\notag
\end{equation}
This proves the first assertion of the proposition.
Now let us assume that $u$ is not derived of $v$ for any $v\in
supp(r)$. Let $i\in\{2,\ldots,n\}$. Since $u=u_1\neq u_i$,
Proposition~\ref{prop2.9} gives the following implications:
\begin{equation}
u\in supp(\psi(u_i))\Rightarrow u\in supp(\psi(u_i)-u_i)\Rightarrow
\text{$u$ is derived of $u_i$}\notag
\end{equation}
By assumption on $u$, this implies that $u^*(\psi(u_i))=0$ for any
$i\geqslant 2$. Using $(i)$, we deduce the announced conclusion:
$u^*(\psi(r))=t_1=u^*(r)\neq 0$
\hfill$\square$\\

\subsection{Ordering products of transvections}

In Proposition~\ref{prop2.11} we have seen that $\psi(\alpha)$ may be
computed easily when $\psi\in \c T$ and $\alpha\in Q_1$ provided that
$\psi$ is written as a decreasing product of transvections. The main
result of this
subsection proves that any $\psi\in\c T$ can be uniquely written
that way. Recall that $<$ is an order on the set of non trivial
paths in $Q$ defined in Definition~\ref{def2.5}. The following
notations will be useful.
\begin{definition}
\label{def2.14}
Let $(\alpha,u)$ be a bypass. We set $\c T_{<(\alpha,u)}$ and $\c
T_{\leqslant (\alpha,u)}$ to be the subgroups of $\c T$ generated by
the following sets of transvections:
\begin{equation}
\begin{array}{l}
\{\varphi_{\beta,v,\tau}\ |\ (\beta,v)<(\alpha,u)\ \text{and}\ \tau\in
k\}\ \ for\ \c T_{<(\alpha,u)}\\
\{\varphi_{\beta,v,\tau}\ |\ (\beta,v)\leqslant (\alpha,u)\ \text{and}\ \tau\in
k\}\ \ for\ \c T_{\leqslant(\alpha,u)}
\end{array}\notag
\end{equation}
Also, we define $\c T_{(\alpha,u)}$ to be the following subgroup of
$\c T$:
\begin{equation}
\c T_{(\alpha,u)}=\{\varphi_{\alpha,u,\tau}\ |\ \tau\in k\}\notag
\end{equation}
\end{definition}
\begin{rem}
\label{rem2.15}
\begin{enumerate}
\item[.] $\c T_{(\alpha,u)}$ is indeed a subgroup of $\c T$ because
  $\varphi_{\alpha,u,\tau}\varphi_{\alpha,u,\tau'}=\varphi_{\alpha,u,\tau+\tau'}$ 
for any $\tau,\tau'\in k$. Actually, the following mapping is an
isomorphism of abelian groups:
\begin{equation}
\begin{array}{rcl}
k&\longrightarrow & \c T_{(\alpha,u)}\\
\tau&\longmapsto & \varphi_{\alpha,u,\tau}
\end{array}\notag
\end{equation}
\item[.] $\c T_{\leqslant(\alpha,u)}$ is generated by $\c
  T_{<(\alpha,u)}\cup \c T_{(\alpha,u)}$.
\item[.] If $(\alpha,u)<(\beta,v)$, then $\c
  T_{\leqslant(\alpha,u)}\subseteq \c T_{\leqslant(\beta,v)}$ and
  $T_{<(\alpha,u)}\subseteq T_{<(\beta,v)}$.
\item[.] $\c T=\bigcup\limits_{(\alpha,u)}\c T_{\leqslant (\alpha,u)}$
  and if $(\alpha_m,u_m)$ is the greatest bypass in $Q$, then $\c T=\c
  T_{\leqslant {(\alpha_m,u_m)}}$ (recall that $Q$ has finitely many
  bypasses because it has no oriented cycle).
\end{enumerate}
\end{rem}

The following lemma proves that any $\psi\in \c T$ is a decreasing
product of transvections.
\begin{lem}
\label{lem2.16}
\begin{enumerate}
\item[.] $\c T_{<(\alpha,u)}$ is a normal subgroup of $\c T_{\leqslant
    (\alpha,u)}$, for any bypass $(\alpha,u)$.
\item[.] Let $(a_1,v_1)<\ldots<(a_N,v_N)$ be the (finite) increasing
  sequence of all the bypasses in $Q$. Then:
\begin{enumerate}
\item[-] $\c T_{ <(a_i,v_i)}=\c T_{\leqslant(a_{i-1},v_{i-1})}$ if
      $i\geqslant 1$,
\item[-] $\c T_{<(a_1,v_1)}=1$,
\item[-] $\c T_{\leqslant (a_i,v_i)}=\c T_{(a_i,v_i)}\c
  T_{(a_{i-1},v_{i-1})}\ldots\c T_{(a_1,v_1)}$.
\end{enumerate}
\end{enumerate}
\end{lem}
\noindent{\textbf{Proof:}} Thanks to Remark~\ref{rem2.15}, it is
sufficient  to prove that if $\tau,\nu\in k$ and if $(\beta,v),(\alpha,u)$
are bypasses such that $(\beta,v)<(\alpha,u)$, then:
\begin{equation}
\varphi_{\beta,v,\nu}\varphi_{\alpha,u,\tau}\in
\varphi_{\alpha,u,\tau}\c T_{<(\alpha,u)}\tag{$\star$}
\end{equation}
There are two situations wether $(\alpha,u,\beta,v)$ is a double
bypass or not. 
If $(\alpha,u,\beta,v)$ is a double bypass, then 
Section~$1$ gives:
\begin{equation}
\varphi_{\beta,v,\nu}\varphi_{\alpha,u,\tau}=
\varphi_{\alpha,u,\tau}\varphi_{\alpha,w,\tau\nu}\varphi_{\beta,v,\nu}\notag
\end{equation}
where $w$ is the path obtained from $u$ after replacing $\beta$ by
$v$. Moreover, Lemma~\ref{lem2.7} implies that
$(\beta,v)<(\alpha,w)<(\alpha,u)$. Therefore, $(\star)$ is satisfied when
$(\alpha,u,\beta,v)$ is a double bypass.
If $(\alpha,u,\beta,v)$ is not a double bypass, then
Section~$1$ gives (notice that thanks to Lemma~\ref{lem2.7} and to
the inequality $(\beta,v)<(\alpha,u)$ we know that $(\beta,v,\alpha,u)$ is not a
double bypass):
\begin{equation}
\varphi_{\alpha,u,\tau}\varphi_{\beta,v,\nu}=\varphi_{\beta,v,\nu}\varphi_{\alpha,u,\tau}\notag
\end{equation} 
So $(\star)$ is also satisfied when $(\alpha,u,\beta,v)$ is not a double
bypass.\hfill$\square$\\

Using the preceding lemma and Proposition~\ref{prop2.11}, it is now
possible to prove that any $\psi\in \c T$ is
uniquely  a decreasing product of transvections.
\begin{prop}
\label{prop2.17}
Let $(\alpha,u)$ be a bypass and let $\psi\in\c
T_{\leqslant(\alpha,u)}$. Then, there exist a non negative integer $n$,
a sequence of bypasses $(\alpha_1,u_1),\ldots,(\alpha_n,u_n)$ and non
zero scalars $\tau_1,\ldots,\tau_n\in k^*$ verifying:
\begin{enumerate}
\item[(i)]
  $\psi=\varphi_{\alpha_n,u_n,\tau_n}\ldots\varphi_{\alpha_1,u_1,\tau_1}$,
\item[(ii)] $(\alpha_1,u_1)<\ldots<(\alpha_n,u_n)\leqslant (\alpha,u)$.
\end{enumerate}
Moreover, the integer $n$ and the sequence
$(\alpha_1,u_1,\tau_1),\ldots,(\alpha_n,u_n,\tau_n)$ are unique for
these properties.
\end{prop}
\noindent{\textbf{Proof:}} 
The existence is given by Lemma~\ref{lem2.16}. So it suffices to
characterise the triples $(\alpha_i,u_i,\tau_i)$ using $\psi$ only. Let
$A,B$ and $T$ be the following sets:
\begin{equation}
\begin{array}{l}
A:=\{\alpha\in Q_1\ |\ \psi(\alpha)\neq \alpha\}\\
B:=\{(\alpha,u)\ |\ \text{$(\alpha,u)$ is a bypass, $\alpha\in A$ and
  $u\in supp(\psi(\alpha))$}\}\\
T:=\{(\alpha,u,\tau)\ |\ (\alpha,u)\in B\ \text{and}\
\tau=u^*(\psi(\alpha))\}
\end{array}\notag
\end{equation}
Notice that the definition of $A,B,T$ depend on $\psi$ only (and not
on the triples $(\alpha_i,u_i,\tau_i)$). Let $\beta\in Q_1$. Then
Proposition~\ref{prop2.11} gives a normal form:
\begin{equation}
\psi(\beta)=\beta+\sum\limits_{i\ \text{such that}\
  \beta=\alpha_i}\tau_iu_i\notag
\end{equation}
By definition of a normal form and because of ($i$) and ($ii$), the
following equalities hold:
\begin{equation}
\begin{array}{l}
A=\{\alpha_1,\ldots,\alpha_n\}\\
B=\{(\alpha_1,u_1),\ldots,(\alpha_n,u_n)\}\\
T=\{(\alpha_1,u_1,\tau_1),\ldots,(\alpha_n,u_n,\tau_1)\}
\end{array}\notag
\end{equation}
This proves that $n$ and
$(\alpha_1,u_1,\tau_1),\ldots,(\alpha_n,u_n,\tau_n)$ are uniquely
determined by the sets $A,B,T$ (which depend on $\psi$ only)  and by
the total order $<$.\hfill$\square$\\

\section{Comparison of the presentations of a monomial  algebra}
Let $A=kQ/I_0$ with $I_0$ a monomial admissible ideal of $kQ$ and let
$kQ/I\simeq A$ be an admissible presentation of $A$. Thanks to
Proposition~\ref{prop1.1}, there exists $\psi$ a
product of transvections and of a dilatation such that
$\psi(I_0)=I$. The aim of this section is to exhibit $\psi_I$ the
``simplest'' possible among all the $\psi$'s verifying
$\psi(I_0)=I$. It will appear that $\psi_I$ verifies a property which
makes it unique. The construction of $\psi_I$ will use specific
properties of the Gröbner basis of $I$, due to the fact that $I_0$ is
monomial. So, throughout the section, $<$ will denote a total order 
on the set of non trivial paths in $Q$, as in
Definition~\ref{def2.5}. Before studying the Gröbner basis of $I$, it
is useful to give some properties on
the automorphisms $\psi\in AUt_0(kQ)$ verifying $\psi(I_0)=I_0$.
\begin{lem}
\label{lem3.1}
Let $D\in \c D$ be a dilatation. Then $D(I_0)=I_0$. As a consequence,
if $kQ/I\simeq A$ is an admissible presentation, then there exists
$\psi\in\c T$ such that $\psi(I_0)=I$.
\end{lem}
\noindent{\textbf{Proof:}} The first assertion is due to the fact that
$D(u)\in k^*u$ for any path $u$ and to the fact that $I_0$ is
monomial. The second one is a consequence of the first one and of
Proposition~\ref{prop1.1}.\hfill$\square$\\
\begin{lem}
\label{lem3.2}
Let $(\alpha,u)$ be a bypass in $Q$. Then exactly one of the two
following assertions is satisfied:
\begin{enumerate}
\item[.] $\varphi_{\alpha,u,\tau}(I_0)=I_0$ for any $\tau\in k$.
\item[.] $\varphi_{\alpha,u,\tau}(I_0)\neq I_0$ for any $\tau\in
  k^*$.
\end{enumerate} 
\end{lem}
\noindent{\textbf{Proof:}} Assume that $\tau\in k^*$ verifies
$\varphi_{\alpha,u,\tau}(I_0)=I_0$ and let $\mu\in k$. Let $v\in I_0$
be a path. If $\alpha$ does not appear in $v$, then
$\varphi_{\alpha,u,\nu}(v)=v\in I_0$. Assume that $\alpha$ appears
in $v$, i.e. $v=v_2\alpha v_1$ with $v_1,v_2$ paths in which $\alpha$
does not appear (because $Q$ has no oriented cycle). Therefore,
$\varphi_{\alpha,u,\tau}(v)=v+\tau v_2uv_1\in I_0$.
Thus, $v_2uv_1\in
I_0$. This implies that $\varphi_{\alpha,u,\nu}(v)=v+\nu v_2uv_1\in
I_0$. Since $I_0$ is monomial, 
$\varphi_{\alpha,u,\nu}(I_0)=I_0$.\hfill$\square$\\
\begin{lem}
\label{lem3.3}
Let $(\alpha_1,u_1)<\ldots<(\alpha_n,u_n)$ be an increasing sequence
of bypasses, let $\tau_1,\ldots,\tau_n\in k^*$ and set
$\psi=\varphi_{\alpha_n,u_n,\tau_n}\ldots\varphi_{\alpha_1,u_1,\tau_1}$.
Then:
\begin{equation}
\psi(I_0)=I_0\ \ \Leftrightarrow\ \  \varphi_{\alpha_i,u_i,\tau_i}(I_0)=I_0\
\text{for any $i$}\notag
\end{equation}
\end{lem}
\noindent{\textbf{Proof:}} Assume that $\psi(I_0)=I_0$. Let
$i\in\{1,\ldots,n\}$, let $u=a_r\ldots a_1\in I_0$ be a path (with
$a_i\in Q_1$) and fix $i\in \{1,\ldots,n\}$.  If $a_j\neq \alpha_i$ for any
 $j\in\{1,\ldots,r\}$ then
$\varphi_{\alpha_i,u_i,\tau_i}(u)=u\in I_0$. Now assume that there
exists $j\in\{1,\ldots,r\}$ such that $a_j=\alpha_i$ ($j$ is
necessarily unique because $Q$ has no oriented cycle). Therefore:
\begin{equation}
\varphi_{\alpha_i,u_i,\tau_i}(u)=u+\tau_ia_r\ldots
a_{j+1}u_ia_{j-1}\ldots a_1\tag{$i$}
\end{equation}
On the other hand, Proposition~\ref{prop2.9} and
Proposition~\ref{prop2.11} imply that $a_r\ldots
a_{j+1}u_ia_{j-1}\ldots a_1\in supp(\psi(u))$. Thus (recall
that $\psi(u)\in I_0$ and that $I_0$ is monomial):
\begin{equation}
a_r\ldots a_{j+1}u_ia_{j-1}\ldots a_1\in I_0\tag{$ii$}
\end{equation}
From $(i)$ and $(ii)$ we deduce that
$\varphi_{\alpha_i,u_i,\tau_i}(u)\in I_0$ for any path $u\in I_0$. So
$\varphi_{\alpha_i,u_i,\tau_i}(I_0)=I_0$ for any $i$. The remaining
implication is immediate.\hfill$\square$\\

\begin{rem}
\label{rem3.4}
The three preceding lemmas imply that the group $Aut_0(kQ,I_0)$
defined as follows:
\begin{equation}
Aut_0(kQ,I_0):=\{\psi\in Aut(kQ)\ |\ \psi(x)=x\ for\ any\ x\in
  Q_0,\ and\ \psi(I_0)=I_0\}\notag
\end{equation}
is generated by the dilatations and by all the transvections preserving
$I_0$:
\begin{equation}
Aut_0(kQ,I_0)=<\c D\cup \{\varphi\ |\ \varphi\ is\ a\ transvection\
such\ that\ \varphi(I_0)=I_0\}>\notag
\end{equation}
\end{rem}

\begin{ex}
  The preceding remark does not hold for any ideal $I$, even if $kQ/I$
  is monomial. For example, let $Q$ be the quiver:
  
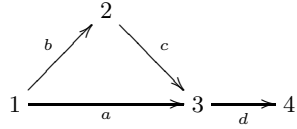
\begin{equation}
    \xymatrix{
      &2\ar@{->}[rd]^c &&\\
      1\ar@{->}[ru]^b \ar@{->}[rr]_a & & 3\ar@{->}[r]_d &4
    }\notag
  \end{equation}
and let $I=<da-dcb>$. Notice that $kQ/I\simeq kQ/I_0$ where
$I_0=\varphi_{a,cb,1}(I)=<da>$. On the other hand:
\begin{enumerate}
\item $Id=\varphi_{a,cb,0}$ is the only transvection lying in
  $Aut_0(kQ,I)$,
\item for $t\in k\backslash\{0,1\}$, the dilatation $D_t$ such that $D_t(a)=ta$ and $D_t(x)=x$ for any
  other arrow  $x$ does not belong to $Aut_0(kQ,I)$,
\item $D_t\varphi_{a,cb,t}\in Aut_0(kQ,I)$ for any $t\in k^*$.
\end{enumerate}
So $Aut_0(kQ,I)$ is not generated by $\c D$ and by the transvections
it contains.
\end{ex}
The following proposition gives the announced properties on the
Gröbner bases of the admissible
ideals $I$ of $kQ$ such that $kQ/I\simeq A$. Recall that for such an
$I$, there exists $\psi\in \c T$ such that $\psi(I_0)=I$ (see
 Lemma~\ref{lem3.1}).
\begin{prop}
\label{prop3.5}
Let $\psi\in\c T$ and let let $I=\psi(I_0)$. Let $B_0$ (resp. $B$) be
the Groebner basis of $I_0$ (resp. of $I$). Then $B_0$ is made of all
the paths in $Q$ which belong to $I_0$. Moreover, the mapping:
\begin{equation}
\begin{array}{rcl}
B&\longrightarrow & B_0\\
r&\longmapsto & max(supp(r))
\end{array}\tag{$\star$}
\end{equation}
is well defined and bijective. For $u\in B_0$, let $r_u\in B$
be the inverse image of $u$ under $(\star)$. Then $supp(r_u-u)$ is a
set of paths derived of $u$.
\end{prop}
\noindent{\textbf{Proof:}} Let $u_1<\ldots<u_n$ be the increasing
sequence of all the non trivial paths in $Q$. Let $(r_1,\ldots,r_d)$
be the Gröbner basis of $I$ and for each $j\in\{1,\ldots,d\}$, let
$i_j\in\{1,\ldots,n\}$ be such that:
\begin{equation}
r_j\in u_{i_j}+Span(u_l\ ;\ l<i_j)\notag
\end{equation}
Since $I_0$ is monomial, $B_0$ is made of all the paths in $Q$
belonging to $I_0$. 

Let $j\in \{1,\ldots,d\}$. Since $u_{i_j}=max(supp(r_j))$, the path
$u_{i_j}$ is not derived of $u$ for any $u\in supp(r_j)$ (thanks to
Lemma~\ref{lem2.7}). So Proposition~\ref{prop2.13} implies that
$u_{i_j}\in supp(\psi^{-1}(r_j))\in I_0$. Because $I_0$ is monomial,
this proves that $u_{i_j}\in I_0$. Therefore, the mapping $(\star)$ is
well defined. It is also one-to-one because of the definition of the
Groebner basis of $I$.
 Let $u\in
B_0$. Proposition~\ref{prop2.9} implies that
$u=max(supp(\psi(u))$. Since $\psi(u)\in I$, there
exists $j\in\{1,\ldots,d\}$ such that $u=u_{i_j}=max(supp(r_j))$. This
proves that $(\star)$ is onto and therefore bijective.

It remains to prove the last assertion of the proposition. This will
be done by proving by induction on $j\in\{1,\ldots,d\}$ that the following
assertion is true:
\begin{equation}
  H_j:"\text{$supp(r_j-u_{i_j})$ is a set of paths derived of
  $u_{i_j}$}"\notag
\end{equation}
Remark that Proposition~\ref{prop2.9} implies that for any $j$:
\begin{equation}
u_{i_j}=max(supp(\psi(u_{i_j})))\ \ \text{and}\ 
u_{i_j}^*(\psi(u_{i_j}))=1\tag{$i$}
\end{equation}
Moreover, $\psi(u_{i_j})\in I$ because $(\star)$ is well defined and
 because $\psi(I_0)=I$.
Now begins the induction. Both $r_1$ and $\psi(u_{i_1})$ lie in
 $I$. Moreover, $u_{i_1}=max(supp(r_1))$ by definition of $u_{i_1}$
 and $u_{i_1}=max(supp(\psi(u_{i_1})))$ because of
 Proposition~\ref{prop2.9}. So $H_1$ is true.
 Assume that $j\geqslant 2$ and that $H_1,\ldots,H_{j-1}$ are
true. Since $\psi(u_{i_j})\in I$ and because of $(i)$, the following holds:
\begin{equation}
\psi(u_{i_j})=r_j+\sum\limits_{
\begin{array}{c}
j'<j,\\
u_{i_{j'}}\in supp(\psi(u_{i_j}))
\end{array}
}
u_{i_{j'}}^*(\psi(u_{i_j}))r_{j'}\notag
\end{equation}
So:
\begin{equation}
r_j-u_{i_j}=\psi(u_{i_j})-u_{i_j}-\sum\limits_{
\begin{array}{c}
j'<j,\\
u_{i_{j'}}\in supp(\psi(u_{i_j}))
\end{array}
}
u_{i_{j'}}^*(\psi(u_{i_j}))\left[(r_{j'}-u_{i_{j'}})+u_{i_j'}\right]
\tag{$ii$}
\end{equation}
Notice that in the above equality:
\begin{enumerate}
\item[$(iii)$] $supp(\psi(u_{i_j})-u_{i_j})$ is a set of paths derived of
  $u_{i_j}$ (thanks to Proposition~\ref{prop2.9}),
\item[$(iv)$] if $j'<j$ verifies $u_{i_{j'}}\in supp(\psi(u_{i_j}))$, then:
\begin{enumerate}
\item[$(v)$] $u_{i_{j'}}$ is derived of $u_{i_j}$ (see $(iii)$ above),
\item[$(vi)$] $supp(r_{j'}-u_{i_{j'}})$ is a set of paths derived of
  $u_{i_{j'}}$ (because $H_{j'}$ is true) and therefore derived of
  $u_{i_j}$ (thanks to $(v)$ and to Lemma~\ref{lem2.3}).
\end{enumerate}
\end{enumerate}
The points $(ii)-(vi)$ prove that $H_j$ is true. Hence,
$H_j$ is true for any $j\in\{1,\ldots,d\}$. This finishes the proof of
the proposition.\hfill$\square$\\

Now it is possible to define precisely the automorphism $\psi_I$ mentionned at the
beginning of the section.
\begin{prop}
\label{prop3.6}
Let $kQ/I\simeq A$ be an admissible presentation. Then there exists a
unique $\psi_I\in\c T$ verifying the following conditions:
\begin{enumerate}
\item[$1)$] $\psi_I(I_0)=I$,
\item[$2)$] if $(\alpha,u)$ is a bypass such that $u\in
  supp(\psi_I(\alpha))$ then $\varphi_{\alpha,u,\tau}(I_0)\neq I_0$
  for any $\tau\in k^*$ (see Lemma~\ref{lem3.2}).
\end{enumerate}
\end{prop}
\noindent{\textbf{Proof:}} $\bullet$ First, the existence of
$\psi_I$. Thanks to Lemma~\ref{lem3.1}, there exists
$\psi\in\c T$ verifying $1)$. Set:
\begin{equation}
\c A := \{\psi\in\c T\ |\ \psi(I_0)=I\}\notag
\end{equation}
and assume that for any $\psi\in\c A$, the condition $2)$ is not
verified. So, for any $\psi\in\c A$, there is a finite (recall that $Q$
has no oriented cycle) and non empty
set of bypasses (see Lemma~\ref{lem3.2}):
\begin{equation}
B_{\psi} = \left\{
(\alpha,u)\ \left|\
\begin{array}{l}
\text{$(\alpha,u)$ is a bypass}\\
u\in supp(\psi(\alpha))\\
\varphi_{\alpha,u,\tau}(I_0)=I_0\ \text{for any $\tau\in k$}
\end{array}\right\}\right.\notag
\end{equation}
For each $\psi\in\c A$, let $(\alpha_{\psi},u_{\psi})=max\ B_{\psi}$ and
let $\psi\in A$ be such that:
\begin{equation}
(\alpha_{\psi},u_{\psi})=min\ \{(\alpha_{\psi'},u_{\psi'})\ |\ \psi'\in\c
A\}\notag
\end{equation}
For simplicity, set $(\alpha,u):=(\alpha_{\psi},u_{\psi})$, $\tau:=u^*(\psi(\alpha))$
$\psi':=\psi\varphi_{\alpha,u,-\tau}$.
Notice that $\psi'\in\c A$ because $(\alpha,u)\in B_{\psi}$. In order to
get a contradiction, let us prove that
$(\alpha_{\psi'},u_{\psi'})<(\alpha,u)$. To do this, let us
prove first that $(\alpha,u)\not\in B_{\psi'}$. Thanks to
Proposition~\ref{prop2.17}, the following equality holds:
\begin{equation}
\psi=\varphi_{\alpha_n,u_n,\tau_n}\ldots\varphi_{\alpha_1,u_1,\tau_1}\notag
\end{equation}
where $(\alpha_1,u_1)<\ldots<(\alpha_n,u_n)$ and where
$\tau_1,\ldots,\tau_n\in k^*$. On the other hand, since
$u^*(\psi(\alpha))=\tau\neq 0$, 
Proposition~\ref{prop2.11} gives:
\begin{equation}
(\exists ! i\in\{1,\ldots,n\})\ \
(\alpha_i,u_i,\tau_i)=(\alpha,u,\tau)\notag
\end{equation}
Let us set:
\begin{equation}
\psi_1:=\varphi_{\alpha_{i-1},u_{i-1},\tau_{i-1}}\ldots\varphi_{\alpha_1,u_1,\tau_1}
\in \c T_{<(\alpha,u)}\notag
\end{equation}
Hence, the following equality holds:
\begin{equation}
\psi'=\varphi_{\alpha_n,u_n,\tau_n}\ldots\varphi_{\alpha_{i+1},u_{i+1},\tau_{i+1}}
\varphi_{\alpha,u,\tau}\psi_1\varphi_{\alpha,u,\tau}^{-1} \notag
\end{equation}
Since $\psi_1\in\c T_{<(\alpha,u)}$, Lemma~\ref{lem2.16} implies that
$\varphi_{\alpha,u,\tau}\psi_1\varphi_{\alpha,u,\tau}^{-1}\in \c
T_{<(\alpha,u)}$. Therefore, Proposition~\ref{prop2.17} gives the
equality:
\begin{equation}
\varphi_{\alpha,u,\tau}\psi_1\varphi_{\alpha,u,\tau}^{-1} =
\varphi_{\beta_m,v_m,\nu_m}\ldots\varphi_{\beta_1,v_1,\nu_1}\notag
\end{equation}
where $(\beta_1,v_1)<\ldots<(\beta_m,v_m)<(\alpha,u)$ and
$\nu_1,\ldots,\nu_m\in k^*$. As a consequence:
\begin{equation}
\psi'=\varphi_{\alpha_n,u_n,\tau_n}\ldots\varphi_{\alpha_{i+1},u_{i+1},\tau_{i+1}}
\varphi_{\beta_m,v_m,\nu_m}\ldots\varphi_{\beta_1,v_1,\nu_1}\notag
\end{equation}
where
$(\beta_1,v_1)<\ldots<(\beta_m,v_m)<(\alpha,u)<(\alpha_{i+1},u_{i+1})<\ldots<(\alpha_n,u_n)$
and where $\tau_{i+1},\ldots,\tau_n,\nu_1,\ldots,\nu_m\in
k^*$. In particular, Proposition~\ref{prop2.11} implies that $u\not\in
supp(\psi'(\alpha))$. Therefore,
$(\alpha,u)\not\in B_{\psi'}$ and in particular,
$(\alpha,u)\neq (\alpha_{\psi'},u_{\psi'})=max\ B_{\psi'}$.
Thus, in order to prove that $(\alpha_{\psi'},u_{\psi'})<
(\alpha,u)$, it suffices to pove that the following
implication holds for any bypass $(\beta,v)$:
\begin{equation}
v\in supp(\psi'(\beta))\ \text{and}\ (\alpha,u)<(\beta,v)\ \Rightarrow
\varphi_{\beta,v,t}(I_0)\neq I_0\ \text{for any $\tau\in k^*$}\tag{$i$}
\end{equation}
Let $(\beta,v)$ be a bypass such that $v\in supp(\psi'(\beta))$ and
such that $(\alpha,u)<(\beta,v)$. Since
$\psi'=\psi\varphi_{\alpha,u,-\tau}$, the following holds:
\begin{equation}
\psi'(\beta)=\left\{
\begin{array}{ll}
\psi(\beta)& \text{if $\beta\neq \alpha$}\\
\psi(\beta)-\tau\psi(u)&\text{if $\beta=\alpha$}
\end{array}\right.
\notag
\end{equation}
Therefore, $v\in supp(\psi'(\beta))\subseteq supp(\psi(\beta))\cup
supp(\psi(u))$. Remark that if $v\in supp(\psi(u))\backslash supp(\psi(\beta))$, then
$\alpha=\beta$ and Proposition~\ref{prop2.9} implies that $v$ is derived of $u$ (we have
$u\neq v$ because $\beta=\alpha$ and $(\alpha,u)<(\beta,v)$) and
therefore $(\alpha,u)>(\alpha,v)=(\beta,v)$ whereas we assumed that
$(\alpha,u)<(\beta,v)$. This proves that $v\in
supp(\psi(\beta))$. Since
$(\beta,v)>(\alpha,u)=(\alpha_{\psi},u_{\psi})= max\ B_{\psi}$ we deduce that
$\varphi_{\beta,v,\tau}(I_0)=I_0$ for any $\tau\in k$. This proves
that the implication $(i)$ is satisfied. Thus:
\begin{equation}
(\alpha_{\psi'},u_{\psi'})< (\alpha,u)=(\alpha_{\psi},u_{\psi})\notag
\end{equation}
This contradicts the minimality of $(\alpha_{\psi},u_{\psi})$ and
proves the existence of $\psi$.

$\bullet$ It remains to prove the uniqueness of $\psi_I$. Assume
that $\psi,\psi'\in\c T$ verify the conditions $1)$ and $2)$. In order
to prove that $\psi=\psi'$, it is sufficient to prove that
$\theta^*(\psi(\alpha))=\theta^*(\psi'(\alpha))$ for any bypass
$(\alpha,\theta)$. Let $\alpha\in Q_1$ and assume that there exists a
minimal path $\theta$ such that $(\alpha,\theta)$ is bypass and such that
$\theta^*(\psi(\alpha))\neq\theta^*(\psi'(\alpha))$. We may assume
that $\theta^*(\psi(\alpha))\neq 0$, i.e. $\theta\in
supp(\psi(\alpha))$. Since $\psi$ verifies $2)$, we deduce that there
exist paths $u$ and $v$ such that:
\begin{equation}
u\in I_0,\ v\not\in I_0\ \text{and}\
\varphi_{\alpha,\theta,1}(u)=u+v\not\in I_0\notag
\end{equation}
Notice that Proposition~\ref{prop2.9} gives:
\begin{equation}
\left\{
\begin{array}{l}
v^*(\psi(u))=\theta^*(\psi(\alpha))\ \text{and}\ u^*(\psi(u))=1\\
v^*(\psi'(u))=\theta^*(\psi'(\alpha))\ \text{and}\ u^*(\psi'(u))=1
\end{array}\right.\tag{$ii$}
\end{equation}
Moreover, $\psi(u),\psi'(u)\in I_0$ because $u\in I$. Therefore,
Proposition~\ref{prop3.5} gives, the same notations concerning the
Groebner bases, we have:
\begin{equation}
\left\{
\begin{array}{l}
\psi(u)=r_u+\sum\limits_{w\in \c A_{\psi}}w^*(\psi(u))r_w\\
\psi'(u)=r_u+\sum\limits_{w\in \c A_{\psi'}}w^*(\psi'(u))r_w
\end{array}
\right.\notag
\end{equation}
where $\c A_{\psi}$ is equal to:
\begin{equation}
\c A_{\psi}:=\{w\in supp(\psi(u))\ |\ w\neq u\ \text{and}\ w\in
I_0\}\notag
\end{equation}
So:
\begin{equation}
\left\{
\begin{array}{l}
v^*(\psi(u))=v^*(r_u)+\sum\limits_{w\in \c A_{\psi}}w^*(\psi(u))v^*(r_w)\\
v^*(\psi'(u))=v^*(r_u)+\sum\limits_{w\in \c A_{\psi'}}w^*(\psi'(u))v^*(r_w)
\end{array}\right.\tag{$iii$}
\end{equation}
Let $w\in \c A_{\psi}$ be such that
$v^*(r_w)\neq 0$, i.e. $v\in supp(r_w)$. Remark
that $v\in supp(r_w-w)$ because $v\not\in I_0$ and
 $w\in I_0$. So:
\begin{enumerate}
\item[.] $v$ is derived of $w$ (thanks to Proposition~\ref{prop3.5}
  and because $v\in supp(r_w-w)$).
\item[.] $v$ is derived of $u$ of order $1$ (because
  $\varphi_{\alpha,\theta,1}(u)=u+v$).
\item[.] $w$ is derived of $u$ (because $w\in\c A_{\psi}$ and thanks
  to Proposition~\ref{prop2.9}).
\end{enumerate}
Using Lemma~\ref{lem2.3}, these three facts imply that:
\begin{equation}
u=u_2\alpha u_1,\ v=u_2\theta u_1\ \text{and}\ w=u_2\theta'u_1\tag{$iv$}
\end{equation}
where $u_1,u_2$ are paths and where $\theta'$ is a path derived of
$\theta$. In particular, $(\alpha,\theta')$ is a bypass such that
$\theta'<\theta$ (see Lemma~\ref{lem2.7}). Therefore, the minimality
of $\theta$ forces $\theta'^*(\psi(\alpha))=\theta'^*(\psi'(\alpha)))$. Moreover,
$(iv)$ and Proposition~\ref{prop2.9} imply that
\begin{equation}
w^*(\psi(u))=\theta'^*(\psi(\alpha))=\theta'^*(\psi'(\alpha))=w^*(\psi'(u))\notag
\end{equation}
Therefore we have proved the following implication:
\begin{equation}
w\in\c A_{\psi}\ \text{and}\ v^*(r_w)\neq 0\ \Rightarrow\ w^*(\psi(u))v^*(r_w)=w^*(\psi'(u))v^*(r_w)\tag{$v$}
\end{equation}
After exchangeing the roles of $\psi$ and $\psi'$, the arguments used
to prove $(v)$ also give the following implication:
\begin{equation}
w\in\c A_{\psi'}\ \text{and}\ v^*(r_w)\neq 0\ \Rightarrow\ w^*(\psi(u))v^*(r_w)=w^*(\psi'(u))v^*(r_w)\tag{$vi$}
\end{equation}
Then, $(iii)$, $(v)$ and $(vi)$ give
$v^*(\psi(u))=v^*(\psi'(u))$. This and $(ii)$ imply that
$\theta^*(\psi(\alpha))=\theta^*(\psi'(\alpha))$, a
contradiction. This proves that $\psi=\psi'$.\hfill$\square$\\

\section{Proof of the main theorem}
Let $A=kQ/I_0$ where $I_0$ is a monomial admissible ideal of $kQ$. The
aim of this section is to prove that the quiver $\g$ of the 
homotopy relations of the admissible presentations of $A$ admits
$\sim_{I_0}$ as unique source. This fact will be used in order to
the existence of the universal cover of $A$. Notice that $\sim_{I_0}$ is
a source of $\g$. Indeed, all minimal relations in $I_0$ are monomial
relations so, for any $\sim_I \in \g_0$ we have
$\gamma\sim_{I_0}\gamma'\Rightarrow \gamma\sim_I\gamma'$. In order to
prove that $\sim_{I_0}$ is the unique source in $\g$ it will be proved
that for any admissible presentation $kQ/I\simeq A$, the decomposition
of $\psi_I$ (given by Proposition~\ref{prop3.6}) into a decreasing
product of transvections (see Proposition~\ref{prop2.17}) defines a
path in $\g$ starting at $\sim_{I_0}$ and ending at $\sim_I$. In this
purpose, the following proposition will be useful.
\begin{prop}
\label{prop4.1}
Let $kQ/I\simeq A$ be an admissible presentation. Then, for any bypass
$(\alpha,u)$:
\begin{equation}
u\in supp(\psi_I(\alpha))\Rightarrow u\sim_I\alpha\notag
\end{equation}
\end{prop}
\noindent{\textbf{Proof:}} For simplicity, set $\psi:=\psi_I$. Thanks
to Proposition~\ref{prop2.17} there is an equality: 
\begin{equation}
\psi=\varphi_{\alpha_n,u_n,\tau_n}\ldots\varphi_{\alpha_1,u_1,\tau_1}\notag
\end{equation}
with $(\alpha_1,u_1)<\ldots<(\alpha_n,u_n)$ and
$\tau_1,\ldots,\tau_n\in k^*$. Thanks Proposition~\ref{prop2.11}
it suffices to prove that $\alpha_i\sim_I u_i$ for any
$i$. This will be done using a decreasing induction on
$m\in\{1,\ldots,n\}$. Let $H_m$ be the assertion:
\begin{equation}
H_m:"\text{$\alpha_i\sim_Iu_i$ for any $i\in\{m,m+1,\ldots,n\}$}"\notag
\end{equation}
$H_{n+1}$ is true because $\{i\ |\ n+1\leqslant i\leqslant n\}$ is
empty. So assume that $H_{m+1}$ is true ($m\in\{1,\ldots,n\}$). In
order to prove that $H_m$ is true, it thus suffices to prove that
$\alpha_m\sim_I u_m$. From
Proposition~\ref{prop2.11}, the path $u_m$ lies in $supp(\psi(\alpha_m))$. Hence,
Proposition~\ref{prop3.6} provides a path $u\in I_0$ such that
$\varphi_{\alpha_m,u_m,1}(u)\not\in
I_0$. Therefore, there exist paths $v_1,v_2$ such that:
\begin{equation}
u=v_2\alpha_mv_1,\ \ v:=v_2u_mv_1\not\in I_0\ \text{and}\ \varphi_{\alpha_m,u_m,1}(u)=u+v\tag{$i$}
\end{equation}
Since $\psi(u)\in I$, there exists a decomposition:
\begin{equation}
\psi(u)=r_1+\ldots+r_N\notag
\end{equation}
where $r_1,\ldots,r_N$ are minimal relations in $I$ with pairwise
disjoint supports. Remark that $u,v\in supp(\psi(u))$ thanks to
Proposition~\ref{prop2.9} and to Proposition~\ref{prop2.11}. Without
loss of generality, it may be
assumed that $v\in supp(r_1)$. Let $i\in\{1,\ldots,N\}$ be such that
$u\in supp(r_i)$. If $i=1$ then $u\sim_I v$ and $(i)$ gives
$\alpha_m\sim_I u_m$. So assume that $i\neq 1$. 
Remark that $\psi^{-1}(r_1)\in I_0$ because $r_1\in I$. Since $I_0$ is
monomial, this also implies 
that $v\not\in supp(\psi^{-1}(r_1))$. And thanks to
Proposition~\ref{prop2.13}, this proves that:
\begin{equation}
\text{there exists $w\in supp(r_1)$ such that $v$ is derived from
  $w$}\tag{$ii$}
\end{equation}
Therefore:
\begin{enumerate}
\item[.] $w$ is derived of $u$ since $w\in supp(r_1)\subseteq
  supp(\psi(u))$ (see Proposition~\ref{prop2.9}, notice that $u\neq w$ because
  $u\not\in supp(r_1)$),
\item[.] $v$ is derived of $w$ (see $(ii)$),
\item[.] $v$ is derived of $u$ of order $1$ (because of $(i)$).
\end{enumerate}
Thanks to Lemma~\ref{lem2.3}, these three points imply that:
\begin{equation}
w=v_2\theta v_1\ \text{and $u_m$ is derived of $\theta$}\tag{$iii$}
\end{equation}
Since $w\in supp(\psi(u))$, the equalities $w=v_2\theta v_1$,
$u=v_2\alpha_m v_1$ and Proposition~\ref{prop2.9} imply that
$\theta\in supp(\psi(\alpha_m))$. Hence, there exists
$j\in\{1,\ldots,n\}$ such that:
\begin{equation}
(\alpha_m,\theta)=(\alpha_j,u_j)\notag
\end{equation}
Since $u_m$ is derived of $\theta$ (see $(iii)$), this last equality gives
$u_j=\theta>u_m$ (see Lemma~\ref{lem2.7}) and therefore $j>m$. On the
other hand,
$H_{m+1}$ is true, so:
\begin{equation}
\alpha_m=\alpha_j\sim_I u_j=\theta\tag{$iv$}
\end{equation}
Finally, $v\sim_I w$, because $r_1$ is a minimal relation in $I$ such
that $v,w\in supp(r_1)$. This together with $(i)$, $(iii)$ and $(iv)$
imply that
$\alpha_m\sim_Iu_m$.
So $H_m$ is true and the induction is finished.\hfill$\square$\\

\begin{rem}
\label{rem4.2}
The preceding proposition proves that $\alpha\sim_I u$ for
any $u\in supp(\psi(\alpha))$. On the other hand, $\sim_{I_0}$ is
weaker than $\sim_I$ (i.e. $\gamma\sim_{I_0}\gamma'\Rightarrow
\gamma\sim_I\gamma'$). These two properties are linked in
general. Indeed, in \cite[Prop. 4.2.35, Prop. 42.36]{lemeur_thesis}
the author has proved that if $I$ is an admissible ideal (non
necessarily monomial) of $kQ$ and if $\psi\in\c T$ is such that
$\alpha\sim_{\psi(I)}u$ for any bypass $(\alpha,u)$ such that $u\in
supp(\psi(\alpha))$, then $\sim_I$ is weaker than $\sim_{\psi(I)}$.
\end{rem}

Now it is possible to provethe existence of a path in $\g$ starting at $\sim_{I_0}$ and ending
at $\sim_I$, whenever $kQ/I\simeq A$.
\begin{prop}
\label{prop4.3}
Let $kQ/I\simeq A$ be an admissible presentation. Let
$(\alpha_1,u_1)<\ldots<(\alpha_n,u_n)$ be the bypasses and
 $\tau_1,\ldots,\tau_n\in k^*$ the scalars such that 
$\psi_I=\varphi_{\alpha_n,u_n,\tau_n}\ldots\varphi_{\alpha_1,u_1,\tau_1}$
(see
Proposition~\ref{prop2.17}).
For each $i\in\{1,\ldots,n\}$, set:
\begin{equation}
I_i:=\varphi_{\alpha_i,u_i,\tau_i}\ldots\varphi_{\alpha_1,u_1,\tau_1}(I_0)\notag
\end{equation}
then, for each $i$, exactly one of the two following situations occurs:
\begin{enumerate}
\item[.] $\sim_{I_{i-1}}$ and $\sim_{I_i}$ coincide,
\item[.] $\varphi_{\alpha_i,u_i,\tau_i}$ induces an arrow
  $\sim_{I_{i-1}}\to \sim_{I_i}$ in $\g$.
\end{enumerate}
In particular, there exists a path in $\g$ starting at $\sim_{I_0}$
and ending at $\sim_{I_n}=\sim_I$.
\end{prop}
\noindent{\textbf{Proof:}} Let $i\in\{1,\ldots,n\}$ and set
$\psi_i:=\varphi_{\alpha_i,u_i,\tau_i}\ldots\varphi_{\alpha_1,u_1,\tau_1}$.
Thus $I_i=\psi_i(I_0)$. Using Proposition~\ref{prop2.11} and
Proposition~\ref{prop3.6} it is easily verified that
$\psi_i=\psi_{I_i}$.
Therefore, Proposition~\ref{prop4.1} applied to $I_i$ gives
$\alpha_i\sim_{I_i}u_i$.
Since $I_i=\varphi_{\alpha_i,u_i,\tau_i}(I_{i-1})$, this proves that
(see Proposition~\ref{prop1.2}) either $\sim_{I_{i-1}}$ and $\sim_{I_i}$ coincide or
$\varphi_{\alpha_i,u_i,\tau_i}$ induces an arrow $\sim_{I_{i-1}}\to
\sim_{I_i}$ in $\g$. Thus, the vertices
$\sim_{I_0},\sim_{I_1},\ldots,\sim_{I_n}=\sim_I$ of $\g$
are the vertices of a path in $\g$ (maybe with repetitions) starting
at $\sim_{I_0}$ and ending at $\sim_I$.\hfill$\square$\\

The preceding proposition and the fact that $\g$ has no oriented
cycle gives immediately the following corollary
which was proved by the author in \cite{lemeur2} in the case of
algebras without double bypass over an algebraically
closed field of characteristic zero.
\begin{cor}
\label{cor4.4}
Let $Q$ be a quiver without oriented cycle and without multiple
arrows. Let $I_0$ be an admissible and monomial ideal of $kQ$ and let
$A=kQ/I_0$. Then the quiver $\g$ of the homotopy relations of the
admissible presentations of $A$ admits $\sim_{I_0}$ as unique source.
\end{cor}

The following example shows that the preceding corollary does not hold
if $Q$ has multiple arrows.
\begin{ex}
  Let $A=kQ/I_0$ where $Q$ is the quiver
  $\xymatrix{1\ar@/^/@{->}[r]^a\ar@/_/@{->}[r]_b&2\ar@{->}[r]^c&3}$
  and $I_0=<ca>$. Then $\g$ is equal to:
  
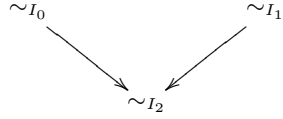
\begin{equation}
    \xymatrix{
      \sim_{I_0} \ar@{->}[rd] &&\sim_{I_1}\ar@{->}[ld]\\
      &\sim_{I_2}&
    }\notag
  \end{equation}
where $I_1=<cb>$ and $I_2=<ca-cb>$. In particular, $\g$ has two
distinct sources. Notice however, that the mapping $a\mapsto b,\
b\mapsto a,\ c\mapsto c$ defines a group isomorphism
$\pi_1(Q,I_0)\simeq \pi_1(Q,I_1)$. One has $\pi_1(Q,I_0)\simeq
\pi_1(Q,I_1)\simeq\mathbb{Z}$ and $\pi_1(Q,I_2)=1$.
\end{ex}
 
Proposition~\ref{prop4.3} also allows one to prove
Theorem~\ref{Thm1}. It extends \cite[Thm. 2]{lemeur2} to monomial
triangular algebras 
without multiple arrows. Notice that Theorem~\ref{Thm1} makes
no assumption on the characteristic of $k$. Also recall that
$\pi_1(Q,I_0)=\pi_1(Q)$.

\noindent{\textbf{Proof of Theorem~\ref{Thm1}:}} The proof is 
identical to the proof of \cite[Thm. 2]{lemeur2} except that one uses
Proposition~\ref{prop4.3} instead of
\cite[Lem. 4.3]{lemeur2}.\hfill$\square$\\

\bibliographystyle{plain}
\bibliography{biblio}
\end{document}